%% file: main.tex
\documentclass[11pt,a4paper]{article}
\usepackage[T1]{fontenc}
\input{glyphtounicode}
\usepackage{lmodern,microtype}
\usepackage{dsfont}
\usepackage[margin=28mm]{geometry}
\usepackage{amsmath,amssymb,amsthm,mathtools}
\usepackage{booktabs,array,longtable,enumitem,needspace}
\usepackage{xcolor}
\usepackage{tikz}
\usetikzlibrary{calc,arrows.meta}
\usepackage{float}
\usepackage[colorlinks=true,linkcolor=blue!45!black,citecolor=blue!45!black,urlcolor=blue!45!black]{hyperref}
\hypersetup{pdftitle={Polytopal realization of subword complexes in Kac-Moody case},pdfauthor={Zijun Li}}
\newtheorem{theorem}{Theorem}[section]
\newtheorem{proposition}[theorem]{Proposition}
\newtheorem{conjecture}[theorem]{Open Problem}
\newtheorem{lemma}[theorem]{Lemma}
\newtheorem{corollary}[theorem]{Corollary}
\theoremstyle{definition}
\newtheorem{definition}[theorem]{Definition}
\newtheorem{example}[theorem]{Example}
\theoremstyle{remark}
\newtheorem{remark}[theorem]{Remark}

\newcommand{\R}{\mathbb R}
\newcommand{\SC}{\operatorname{SC}}
\newcommand{\Dem}{\operatorname{Dem}}
\newcommand{\cone}{\operatorname{cone}}
\newcommand{\conv}{\operatorname{conv}}
\newcommand{\Span}{\operatorname{span}}

\newcommand{\Fac}{\operatorname{Fac}}
\newcommand{\Vertices}{\operatorname{Vert}}

\setlist[enumerate]{itemsep=3pt,topsep=5pt}
\numberwithin{equation}{section}
\newcommand{\SCQ}{\operatorname{SC}(Q,\pi)}
\newcommand{\Z}{\mathbb Z}
\newcommand{\C}{\mathbb C}

\newcommand{\Spec}{\operatorname{Spec}}
\newcommand{\wt}{\widetilde}

\title{Polytopal realization of double root free subword complexes}
\author{Zijun Li}
\date{13 November 2026}
\begin{document}

\maketitle
\begin{abstract}
We give a convex geometric sufficient condition for a convex polytope to realize a spherical subword complex. Using this criterion, we first prove that each root independent spherical subword complex is polytopal, for any Coxeter group. Then we prove the polytopality of double root free subword complexes, for any Coxeter groups. To prove this, we construct a polytope in an auxiliary vector space whose image under a natural affine map is the usual brick polytope. In the complex Kac-Moody case, the polytope is integral, with vertices in the character lattice of the torus of unfolded galleries in a thick building. In finite Lie types, in double root free cases, these tori are precisely the corresponding braid varieties.
\end{abstract}

\section{Introduction}

Knutson and Miller \cite{KM} introduced the concept of subword complexes in the study of Schubert polynomials. Given a Coxeter system $(W,S)$, a word $Q$ whose letters are simple reflections indexed by $S$, and an element $\pi\in W$, the subword complex $\SC_W(Q,\pi)$ has facets given by the complements of reduced expressions of $\pi$ in $Q$. If there is no risk of confusion, we simply denote it by $\SCQ$. Knutson and Miller proved that $\SCQ$ is a simplicial sphere if and only if $\pi$ is the Demazure product of $Q$ and raised the following problem: 

\begin{conjecture}\label{conjecture}
Is any spherical subword complex the boundary complex of a simplicial convex polytope?
\end{conjecture}

Many interesting families of polytopes or simplicial complexes can be realized from the point of view of subword complexes. In the finite case, for a Coxeter element $c$, the subword complex $\SC(cw_0(c),w_0)$ associated with the $c$-sorting word for the longest element $w_0$ is the corresponding cluster complex \cite{CLS}. Its polytopal realizations include the generalized associahedra \cite{FZAssoc}. The boundary of certain cyclic polytopes can be also identified with subword complexes as in \cite{CLS}.

Among different polytopal realizations of subword complexes, the construction of the brick polytopes is of particular interest. Pilaud and Stump introduced brick polytopes for spherical subword complexes of finite Coxeter groups in \cite{PSbrick}. Their construction gave the polytopal realization of spherical root independent subword complexes in the finite case. Jahn and Stump further studied brick polyhedra and the chamber map in the finite Coxeter case to investigate the root dependent case in \cite{JS}.

Another approach to polytopality is given by transformations of subword complexes under braid and $0$-Hecke moves. Gorsky described the corresponding operations on subword complexes and their dual polytopes in \cite{GorskyNilHecke} and \cite{GorskyBraid}. Braid moves can be described by edge subdivisions of subword complexes and by truncations of codimension $2$ faces of polytopes.

A weaker problem for Open Problem~\ref{conjecture} is to ask for a simplicial complete fan realization of the subword complex. Reading and Speyer used Cambrian fans, which are normal fans of generalized associahedra \cite{ReadingSpeyerCambrian}, to study $c$-sorting subword complexes, and extended sortable elements and Cambrian structures to infinite Coxeter groups \cite{ReadingSpeyerInfinite}. In affine types, they gave a description of the $g$-vector fan of acyclic cluster algebras \cite{ReadingSpeyerAffine}. Manneville constructed complete simplicial fan realizations of certain $2$-associahedra using subword complexes
\cite{Manneville2Associahedra}.

For a spherical subword complex $\SCQ$, where $Q=q_1\cdots q_m$ and $d=m-\ell(\pi)$, root independence implies that $d$ is at most the rank of the Coxeter group. Recently, Bossinger, Gorsky and Simental \cite{BGS} studied the polytopality problem using torus actions on compactified braid varieties for finite Weyl groups. They showed that the corresponding braid variety is a torus in the double root free case. The moment polytopes of the torus actions arising from the cluster structures in braid varieties give polytopal realizations of these subword complexes. The dimension of a double root free subword complex is bounded above by the number of positive roots, which is much larger than the rank of the Coxeter group; this upper bound is tight.

Motivated by the work \cite{BGS}, the purpose of this paper is to extend the brick polytope approach to double root free spherical subword complexes in the case where $W$ is infinite. We mainly treat the case where $W$ is the Weyl group of a split Kac-Moody group. And we show that all arguments can be applied to an arbitrary Coxeter system using the Soergel real reflection realization in Proposition~\ref{prop:general-realization} and Remark~\ref{rem:general-Coxeter}. We first prove that the usual brick polytope realizes any root independent spherical subword complex of an arbitrary Coxeter type. Then we show that every double root free spherical subword complex admits a polytopal realization (Theorem~\ref{thm:drf-polytope}). More precisely, if $Q=q_1\cdots q_m$ and $d=m-\ell(\pi)$, we can construct a simple $d$-dimensional polytope, which is called the lifted brick polytope, whose vertices are indexed by the facets of $\SC(Q,\pi)$ and whose edges correspond to flips. Using a similar argument as in \cite{JS}, we can show that its polar dual realizes the subword complex. In the Kac-Moody case, such a polytope is an integral polytope.

Our proof is combinatorial and uses a generalization of the convex geometry argument for brick polytopes in \cite{PSbrick,JS}. This is different from the moment map approach in \cite{BGS}. We associate a basis of a vector space with each facet of the subword complex and give a sufficient condition (Proposition~\ref{prop:convexity-criterion}) that these bases induce a polytopal realization of the subword complex. We first show that this condition is satisfied by the root configuration of facets of a root independent spherical subword complex in the real realization space. In the double root free case, the difficulty is that the dimension of the real realization space is in general smaller than the cardinality $d$ of the facet of the subword complex. We construct an auxiliary lattice $M_Q$ of rank $d$ with a basis $(e_i^I)_{i\in I}$ for each facet $I$ in the Kac-Moody case. These bases are given by a pullback via a natural projection
\[
 p:M_Q\longrightarrow\Lambda_Q,\qquad p(e_i^I)=-r(I,i),
\]
where $\Lambda_Q$ is the root configuration lattice and $r(I,-)$ is the root function associated with the facet $I$. Double root freeness implies that this projection is well defined. We then lift the brick vectors to $M_Q\otimes_{\mathbb Z}\mathbb R$ and apply the convexity criterion (Proposition~\ref{prop:convexity-criterion}) to give the desired polytopal realization. For an arbitrary Coxeter system, such a lift works for vector spaces instead of integral lattices.

The auxiliary lattice $M_Q$ also has an algebraic interpretation in the Kac-Moody case. The subword complex $\SCQ$ is associated with the space $\mathcal G(Q)$ of unfolded galleries with a fixed endpoint $C_{\pi}$ in a thick building of the Coxeter system $(W,S)$. When $W$ is the Weyl group of a minimal split Kac-Moody group over $\mathbb C$, we consider the building associated with this Kac-Moody group. The space of unfolded galleries can be described inside a fiber of a Bott-Samelson variety. Such gallery spaces are related to braid varieties and their brick compactifications as considered in \cite{EscobarBrick,CGGLSC,BGS}. In particular, the space $\mathcal G(Q)$ is an algebraic torus in the double root free case. We identify $M_Q$ with the character lattice of the algebraic torus $\mathcal G(Q)$ and show that the projection $p$ agrees with the pullback of characters under the Cartan action. 

Our criterion~\ref{prop:convexity-criterion} is independent of the specific brick vectors used in \cite{JS}. In fact, Proposition~\ref{prop:convexity-criterion} applies to a real
vector space equipped with a basis for each facet and compatible flip relations. The brick polytope realization in the root independent case in \cite{JS} and our lifted realization in the auxiliary space $M_Q$ in the double root free case are both special cases of this construction. We expect that this construction of the polytopal realization can be applied for more general spherical subword complexes.

In Section~\ref{sec:coxeter}, we recall basic properties of subword complexes and their root functions. In Section~\ref{sec:km}, we give the convexity criterion and prove the polytopal realization in the root independent case. In Section~\ref{sec:double}, we prove the polytopal realization in the double root free case. In Section~\ref{sec:complex}, we give the interpretation of $M_Q$ in terms of the character lattice of the torus $\mathcal G(Q)$. Finally, we give some examples of the brick polytopes and lifted brick polytopes for both finite and affine cases in Section~\ref{sec:interleaved-examples}.

\section{Preliminaries on Subword Complexes}
\label{sec:coxeter}

We fix our conventions for Coxeter systems, see \cite{Hum,Kac,Kumar} for references. Let $(W,S)$ be a Coxeter system, with
\[
 W=\langle S\mid s^2=e,\ (st)^{m_{st}}=e
       \text{ if }s\ne t\text{ and }m_{st}<\infty\rangle,
\]
where $m_{st}=m_{ts}\in\{2,3,\ldots\}\cup\{\infty\}$.
For any $w\in W$, the length $\ell(w)$ is the smallest number of simple reflections expressing $w$. Such an expression of this length is called a \textbf{reduced expression}.
We always assume that $S$ is finite.

For a Coxeter system $(W,S)$, a \textbf{real reflection realization} is a vector space $E$ equipped with linearly independent simple roots $\alpha_s\in E$ and linearly independent simple coroots $\alpha_s^\vee\in E^*$ such that $\langle\alpha_s,\alpha_s^{\vee}\rangle=2$ and the action of $W$ on $E$ can be expressed as
\begin{equation}\label{eq:simple-reflection-action}
 s\lambda=\lambda-\langle\lambda,\alpha_s^\vee\rangle\alpha_s.
\end{equation}
Here the brackets denote evaluation $\langle\lambda,h\rangle=h(\lambda) \text{ for all } \lambda\in E,\ h\in E^*$ and $s$ ranges over all simple reflections. Set $V=\bigoplus_{s\in S}\R\alpha_s\subseteq E$.

The action of $W$ on $E^*$ is contragredient, we have
$\langle w\lambda,wh\rangle=\langle\lambda,h\rangle$.
We choose fundamental weights $\omega_t\in E$ satisfying
$\langle\omega_t,\alpha_s^\vee\rangle=\delta_{st}$. In such a real reflection realization, since all coroots are linearly independent, we can always extend $\{\alpha_s^\vee\}_{s\in S}$ to a basis of $E^*$ and take a dual basis to construct the fundamental weights. We use the following realizations $E$ in this paper.

In finite type, we take $E=V=\bigoplus_{s\in S}\R\alpha_s$ with a $W$-invariant positive definite inner product $(\ ,\ )$. The fundamental weights then belong to $V$. In finite Weyl type, we use the usual crystallographic root realization \cite{Kac}.

In Kac-Moody type, we choose a real realization
$(\mathfrak h,(\alpha_s),(\alpha_s^\vee))$ of the generalized Cartan matrix $A=(a_{st})$, with
\[
 E=\mathfrak h^*,\qquad E^*=\mathfrak h,
 \qquad \langle\alpha_t,\alpha_s^\vee\rangle=a_{st}.
\]
The precise construction is given in \cite{Kac}. The fundamental weights need not lie in $V$ in general. 

For an arbitrary Coxeter system, we use a minimal Soergel realization in \cite[Section~2]{SW} and \cite[Chapter~2]{RicheLectures}. Set $n=|S|$ and $C=(c_{st})_{s,t\in S}$, where
\begin{equation}\label{eq:coxeter-realization}
c_{st}=
\begin{cases}
2, & s=t,\\
-2\cos(\pi/m_{st}), & s\ne t,\ m_{st}<\infty,\\
-2, & m_{st}=\infty.
\end{cases}
\end{equation}
In \cite{SW} and \cite[Chapter~2, \S1.2.2]{RicheLectures}, the authors show that there exists a faithful representation $E$ of dimension $2n-\operatorname{rank}C$, with $n$ linearly independent simple roots $\{\alpha_s\}_{s\in S}$ in $E$ and $n$ linearly independent simple coroots $\{\alpha_s^\vee\}_{s\in S}$ in $E^*$. The action of the Coxeter group $W$ on $E$ is just given by Equation~\eqref{eq:simple-reflection-action}. Moreover, we have $\langle\alpha_t,\alpha_s^\vee\rangle=c_{st}$ for all $s,t\in S$. Its restriction to $V=\bigoplus_{s\in S}\mathbb R\alpha_s$ is the standard geometric representation in \cite{Kac}. 

Since the simple coroots are linearly independent, the map
\[
E\longrightarrow\mathbb R^n,\qquad
\lambda\longmapsto
\bigl(\langle\lambda,\alpha_s^\vee\rangle\bigr)_{s\in S}
\]
is surjective. We choose fundamental weights $\omega_t\in E$
with $\langle\omega_t,\alpha_s^\vee\rangle=\delta_{st}$.

Formula~\eqref{eq:simple-reflection-action} then gives
$s\alpha_t=\alpha_t-c_{st}\alpha_s$ and $s\omega_t=\omega_t-\delta_{st}\alpha_s$.

In any of these realizations, set
\[
 \Phi=W\{\alpha_s:s\in S\},\qquad
 \Phi^+=\Phi\cap\sum_{s\in S}\R_{\ge0}\alpha_s,
 \qquad \Phi^-= -\Phi^+.
\]
In Kac-Moody type, $\Phi$ denotes only the real roots.
For $\beta=w\alpha_s$, its coroot and corresponding reflection can be expressed as
\begin{equation}\label{eq:real-root-reflection}
 \beta^\vee=w\alpha_s^\vee,\qquad s_\beta=wsw^{-1},
 \qquad
 s_\beta\lambda
   =\lambda-\langle\lambda,\beta^\vee\rangle\beta.
\end{equation}
A root $\beta$ is positive in $\Phi$ if and only if its coroot $\beta^{\vee}$ is positive in the coroot system.
We have
\begin{equation}\label{eq:root-length}
 \ell(ws)>\ell(w)
 \quad\Longleftrightarrow\quad w\alpha_s\in\Phi^+
 \quad\Longleftrightarrow\quad
 w\alpha_s^\vee\in\Phi^{\vee,+}.
\end{equation}
In finite case, this is the classical result in \cite{Hum}. For this equivalence relation in an arbitrary Coxeter system, see \cite{FU} for details.

For a word $Q=q_1\cdots q_m$ in $S$, its Demazure product is
defined recursively by $\delta_0=e$ and
\begin{equation}\label{eq:demazure-product}
 \delta_k=
 \begin{cases}
  \delta_{k-1}q_k,&\ell(\delta_{k-1}q_k)>\ell(\delta_{k-1}),\\
  \delta_{k-1},&\ell(\delta_{k-1}q_k)<\ell(\delta_{k-1}),
 \end{cases}
 \qquad \Dem(Q)=\delta_m.
\end{equation}
This differs in general from the ordinary product $q_1\cdots q_m$.
For a subword we retain the order of the selected positions.
We use the chamber model of the Coxeter complex, with chambers
$C_w$ indexed by $w\in W$ and type-$s$ adjacency between
$C_w$ and $C_{ws}$; the gallery terminology is recalled in
Section~\ref{sec:complex}.

Let $Q=q_1\cdots q_m$ be a finite word and put $\pi=\Dem(Q)$ and $d=m-\ell(\pi)$. A facet of $\SC(Q,\pi)$ is the complement in $[m]$ of a reduced expression for $\pi$. Its faces are all subsets of facets. We only consider such \textbf{spherical} subword complexes in this paper. By Knutson-Miller \cite[Corollary~3.8]{KM}, this complex is a simplicial sphere of dimension $m-\ell(\pi)-1$. Every facet position is flippable, and the flip graphs of the complex and of its links are connected. Links are again subword complexes, obtained by deleting the positions of the face from $Q$.

For any subset $I\subseteq[m]$, define
\[
 x_k^I=\prod_{\substack{1\le r\le k\\r\notin I}}q_r,
 \qquad x_0^I=e,\qquad
 r(I,k)=x_{k-1}^I\alpha_{q_k},\qquad \sigma_I=x_m^I,
\]
where the factors appear in the product in the same order as in the word $Q$. In the Coxeter complex, label chambers $C_w=wC_e$. For a facet $I$, we write
\[
 R(I)=\bigl(r(I,i)\bigr)_{i\in I}.
\]
The roots at $[m]\backslash I$ form the inversion sequence of a reduced word for $\pi$; they are distinct positive roots in $\Phi^+\cap\pi\Phi^-$. We write $|\beta|$ for the positive one in $\{\beta,-\beta\}$.

\begin{lemma}\label{lem:flip}
Let $I,J$ be adjacent facets, with $I\setminus\{i\}=J\setminus\{j\}$, and put $\beta=r(I,i)$. Then
\begin{equation}\label{eq:flip-root}
 r(I,j)=\begin{cases}\beta,&i<j,\\-\beta,&j<i,\end{cases}
 \qquad
 r(J,k)=\begin{cases}s_\beta r(I,k),&\min(i,j)<k\le\max(i,j),\\
 r(I,k),&\text{otherwise}.
 \end{cases}
\end{equation}
In particular $r(J,j)=-\beta$, and $\beta\in\Phi^+$ if and only if $i<j$. The flip partner $j$ is the unique crossing position with root $|\beta|$.
\end{lemma}
\begin{proof}
 This is the root function flip rule of \cite[Lemma~3.3 and Lemma~3.6]{CLS}.
\end{proof}

\begin{definition}
A subword complex $\SC(Q,\pi)$ is \textbf{spherical} if $\Dem(Q)=\pi$. 
A subword complex is \textbf{root independent} if the root configuration of any facet is linearly independent.
\end{definition}

We can use reflections with respect to root functions $r(I,-)$ to describe products of letters in subwords.
\begin{lemma}\label{lem:insertion}
Let $I$ be a facet, $J\subseteq I$, and $I\setminus J=\{k_1<\cdots<k_r\}$. Then
\[
 \sigma_J=s_{r(I,k_1)}\cdots s_{r(I,k_r)}\sigma_I.
\]
\end{lemma}

\begin{proof}
Write
\[
\sigma_I=A_0A_1\cdots A_r,
\]
where
\[
A_0:=\prod_{\substack{j<k_1\\ j\notin I}} q_j,\qquad
A_a:=\prod_{\substack{k_a<j<k_{a+1}\\ j\notin I}} q_j
\quad (1\le a<r),
\]
and
\[
A_r:=\prod_{\substack{j>k_r\\ j\notin I}} q_j.
\]

For \(1\le a\le r\), set $u_a:=A_0\cdots A_{a-1}$.
Since \(k_a\) belongs to $I$, we have $s_{r(I,k_a)}=u_a q_{k_a}u_a^{-1}$. Using $u_a^{-1}u_{a+1}=A_a$, we have $\sigma_J = s_{r(I,k_1)}\cdots s_{r(I,k_r)}\sigma_I$.
\end{proof}

In this paper, besides the root independent subword complexes, we focus on another family of subword complexes, which is first introduced in \cite{PVS}:
\begin{definition}\label{def:drf}
    A subword complex is \textbf{double root free} if every subword of $Q$ with ordinary product $\pi$ is reduced.
\end{definition}
The following theorem is part of results in \cite[Theorem~5.9]{BGS} for combinatorial properties of double root free subword complexes; it shows that our definition of double root free words is equivalent to that of Pilaud-Stump in \cite{PVS}.
\begin{theorem}\label{thm:double}
For any Coxeter system, the following conditions are equivalent for a spherical subword complex $\SC(Q,\pi)$:
\begin{enumerate}[label=\textup{(\roman*)}]
\item every subword of $Q$ with ordinary product $\pi$ is reduced;
\item no facet contains two roots equal up to sign;
\item no facet contains two roots equal.
\end{enumerate}
\end{theorem}
\begin{proof}
For (i)$\Rightarrow$(ii), suppose that a facet $I$ has two distinct folded positions $a<b$ with $r(I,a)=\pm r(I,b)$. We have $s_{r(I,a)}=s_{r(I,b)}$. Set $J=I\setminus\{a,b\}$. Lemma~\ref{lem:insertion} gives
\[
 \sigma_J=s_{r(I,a)}s_{r(I,b)}\sigma_I=\pi.
\]
The subword $\{1 \cdots m\}\backslash J$ is longer than the subword $\{1 \cdots m\}\backslash I$, so they can not be reduced at the same time, contradicting (i).

For (ii)$\Rightarrow$(i), suppose that a nonreduced subword of $Q$ has product $\pi$. We use the Coxeter deletion condition to remove two letters at the same time repeatedly until a reduced expression of $\pi$. Let $a<b$ be the two positions removed in this last step and $I$ be the complement in $\{1 \cdots m\}$ of the final reduced expression. The two positions $a,b$ are in the facet $I$. By Lemma~\ref{lem:insertion}, We have
\[
 s_{r(I,a)}s_{r(I,b)}\pi=\pi.
\]
Thus the two roots $r(I,a)$ and $r(I,b)$ define the same reflection. We have $r(I,a)=\pm r(I,b)$, contradicting (ii).

(ii)$\Rightarrow$(iii) is trivial. To prove (iii)$\Rightarrow$(ii), suppose that (iii) holds and that a facet $I$ has $r(I,a)=-r(I,b)$ at positions $a<b$. Set $\eta=|r(I,a)|=|r(I,b)|\in\Phi^+$. Since $\SC(Q,\pi)$ is spherical, there exists a unique $k\in \{1 \cdots m\} \backslash I$ such that $r(I,k)=\eta$. Lemma~\ref{lem:flip} implies that the root function at a folded position is a positive root if and only if this folded position appears before its flip partner. Since $r(I,a)=-r(I,b)$, the only possibility is that $a< k < b$ and $r(I,a)=-r(I,b)=\eta$. Now we consider the flip given by $\{a,k\}$ for the facet $I$. Suppose that $I\backslash \{a\}=J\backslash \{k\}$. We have $r(J,k)=r(J,b)=-\eta$ and $k$ and $b$ are both folded positions in the facet $J$. Contradiction with (iii).

In conclusion, we have (i)$\Leftrightarrow$(ii)$\Leftrightarrow$(iii).

\end{proof}
We choose a special facet of a spherical subword complex, which will serve as a reference facet in the construction of the polytope in Section~\ref{sec:km}.
\begin{lemma}\label{lem:chamber-maps-for-e}
   For a spherical subword complex $\SC(Q,\pi)$, there exists a unique facet $I^-$ with $R(I^-)\subseteq\Phi^-$ and a unique facet $I^+$ with $R(I^+)\subseteq\Phi^+$.
\end{lemma}

\begin{proof}
Suppose that $Q=q_1\cdots q_m$ and $\ell(\pi)=l$. We construct a reduced expression of $\pi$ inductively. Set $x_0=e$. Suppose that $x_k=q_{i_1}\cdots q_{i_k}$ is a reduced subword of $Q$ for any $1\leq k \leq l-1$. We define $x_{k+1}=q_{i_1}\cdots q_{i_k}q_{i_{k+1}}$, where $i_{k+1}$ is the smallest index larger than $i_k$ such that $\ell(x_k q_{i_{k+1}})=k+1$. The complement of the set $\{i_1,\cdots i_l\}$ in $\{1,\cdots, m\}$ is a facet of $\SCQ$, we denote it by $I^-$. By the construction, we have $r(I^-,i)\in \Phi^-$ for any $i\in I^-$. 

Suppose that $I$ is another facet with $R(I)\subseteq\Phi^-$. By the property of the inversion set, for $i\in \{1,\cdots , m\}\backslash I$, we have $r(I,i)\in \Phi^+$. Thus we have $\ell(x^I_{k-1}q_k)=\ell(x^I_k)-1$ if and only if $k\in I$. An induction on $k$ gives $I=I^-$.

Let $Q^{\mathrm{rev}}=q_m\cdots q_1$. We have $\Dem(Q^{\mathrm{rev}})=\pi^{-1}$, and $i\mapsto m+1-i$ gives a bijection between $\Fac \SC(Q,\pi)$ and $\Fac \SC(Q^{\mathrm{rev}},\pi^{-1})$. Write $I^{\mathrm{rev}}$ for the facet in $\SC(Q^{\mathrm{rev}},\pi^{-1})$ corresponding to $I$ in $\SCQ$.  Write $\pi=x_{i-1}^Iv$. We have $\pi^{-1}x_{i-1}^I=v^{-1}$. For $i\in I$, we have
\[
 r(I^{\mathrm{rev}},m+1-i)=\pi^{-1}r(I,i),
\]
where the left-hand side is computed in $Q^{\mathrm{rev}}$. 

By Lemma~\ref{lem:flip}, for any $i\in I$, the root $|r(I,i)|$ belongs to $\Phi^+\cap\pi\Phi^-$, since the subword complex is spherical. Thus
\[
 R(I)\subseteq\Phi^+
 \quad\Longleftrightarrow\quad
 R(I^{\mathrm{rev}})\subseteq\Phi^-.
\]
The root configuration in the right hand side in computed in $\SC(Q^{\mathrm{rev}},\pi^{-1})$. The unique facet $(I^{\mathrm{rev}})^-$ for $\SC(Q^{\mathrm{rev}},\pi^{-1})$ such that $R((I^{\mathrm{rev}})^-)\subseteq \Phi^-$ gives the unique facet $I^+$ of $\SC(Q,\pi)$ such that $R(I^+)\subseteq \Phi^+$.
\end{proof}

\section{Convex criterion for polytopal realization}
\label{sec:km}

In this section, let $(W,S)$ be the Weyl group of a generalized Cartan matrix $A$, with the real realization and chosen fundamental weights in Section~\ref{sec:coxeter}. 

Now we fix a spherical subword complex $\SC(Q,\pi)$ of type $(W,S)$. For a facet $I$ of $\SC(Q,\pi)$, define
\begin{equation}\label{eq:km-weights}
w(I,k)=x_{k-1}^I\omega_{q_k} \qquad \text{and} \qquad
B(I)=\sum_{k=1}^m w(I,k).
\end{equation}
Here $w(I,\cdot)$ is called the \textbf{weight function} of the facet $I$ and $B(I)$ is called the \textbf{brick vector} of the facet $I$ in \cite[Definition~4.1]{PSbrick}.
\begin{definition}\label{def:brick polytope}
    The convex polytope $\mathcal B(Q)=\conv\{B(I):I\in\Fac\SC (Q,\pi)\}$ in $E$ is called the \textbf{brick polytope} of $\SC(Q,\pi)$. The point $\mathbf{b}_Q:=\sum_{k=1}^m \omega_{q_k}$ is called the \textbf{base point} of $\mathcal B(Q)$.
\end{definition}

The following lemma is first proven in \cite[Lemma~4.4 and Lemma~4.5]{PSbrick}.
\begin{lemma}\label{lem:weight-flip}
Fix a pair of adjacent facets $I\setminus\{i\}=J\setminus\{j\}$ of $\SC(Q,\pi)$. For any $k\in [m]$, we have 
\begin{equation}\label{eq:km-weight-flip}
w(J,k)-w(I,k)=-a_{I,i,k}r(I,i),
\end{equation}
for some integer $a_{I,i,k}\geq 0$. Moreover, we have 
\begin{equation}\label{eq:km-brick-flip}
B(J)-B(I)=-\lambda_{I,i}r(I,i),
\end{equation}
where $\lambda_{I,i}=\sum_k a_{I,i,k}>0$.
\end{lemma}

\begin{proof}
Suppose first that $i<j$, and put $\beta=r(I,i)$.
For any $k\notin \{i+1,\cdots,j\}$, we have $x^I_{k-1}=x^J_{k-1}$, thus $w(J,k)=w(I,k)$. When $k\in \{i+1,\cdots,j\}$, write $x_{k-1}^I=x_{i-1}^Iv$, and set $a_{I,i,k}=\langle w(I,k),\beta^\vee\rangle$. The word $q_iv$ is reduced, as in the word $v^{-1}q_i$. Thus the coroot $v^{-1}\alpha^\vee_{q_i}$ is positive. Notice that $a_{I,i,k}=\langle w(I,k),\beta^\vee\rangle=\langle\omega_{q_k},v^{-1}\alpha_{q_i}^\vee\rangle$. The pairing between a fundamental weight and a positive coroot is always a natural number. We have
\begin{equation}\label{eq:coroot-coefficient}
\begin{aligned}
w(J,k)-w(I,k) &=s_{\beta}w(I,k)-w(I,k)\\
              &=-\langle w(I,k),\beta^\vee\rangle\beta\\
              &=-a_{I,i,k}\beta.
\end{aligned}
\end{equation}
In particular, the flip rule $r(I,j)=\beta$ implies $a_{I,i,j}=1$.
Thus $\lambda_{I,i}$ is strictly positive. The case $i>j$ can be proved similarly.
\end{proof}

By Lemma~\ref{lem:flip}, the vector space generated by the root configuration of a facet of $\SC(Q,\pi)$ is not changed under the flip. We define the \textbf{root configuration space} and the \textbf{root configuration lattice} as follows, which are independent of the facet $I$:
\[
E_Q=\Span_{\R}R(I),\qquad \Lambda_Q=\Span_{\Z}R(I).
\]

Now we prove the polytopality of a root independent spherical subword complex in the Kac-Moody case. The proof is a generalization of the method applied for the finite case in \cite{JS}.
\begin{theorem}\label{thm:km-ri}
When $\SC(Q,\pi)$ is root independent, the brick polytope $\mathcal B(Q)$ is a simple polytope in the affine subspace
$B(I^+)+E_Q$ of $\mathfrak h^*$. The brick vector gives a bijection between the vertices of $\mathcal B(Q)$ and facets of $\SC(Q,\pi)$. We have
\begin{equation}\label{eq:km-tangent}
\cone_{B(I)}\mathcal B(Q)=\cone\{-r(I,i):i\in I\}.
\end{equation}
The inclusion-reversing bijection
\begin{equation}\label{eq:km-faces}
F\longmapsto P_F=\conv\{B(I):F\subseteq I\in\Fac\SC (Q,\pi)\}
\end{equation}
gives an isomorphism of simplicial complexes between the subword complex $\SC (Q,\pi)$ and the boundary of the polar dual of $\mathcal B(Q)$.
\end{theorem}

We give a more general argument as a generalization of the argument used for brick polytopes in \cite{PSbrick,JS}. We will use this argument again in the double root free case.
We abstract the conditions the brick vectors need to satisfy to give a polytopal realization of the subword complex. Then we construct the desired vectors for the root independent case and double root free case independently.

Suppose that the cardinality of any facet of the spherical subword complex $\SC(Q,\pi)$ is $d$. Without loss of generality, we suppose that $d>0$. Let $V$ be a $d$-dimensional vector space. Suppose that there exists a basis $(u^I_i)_{i\in I}$ for any facet $I$ of $\SC(Q,\pi)$ such that the following condition is satisfied
at any flip $I\setminus\{i\}=J\setminus\{j\}$ and any $k\in I\cap J$:
\begin{equation}\label{eq:abstract-transport}
u_j^J=-u_i^I,\qquad
u_k^J=u_k^I+c_k u_i^I \text{ for some } c_k\in \mathbb{R} .
\end{equation}
We now construct a polytopal realization of $\SC(Q,\pi)$ based on these bases.
We say that a functional $f$ on $V$ is \textbf{regular} if $f(u^I_i)\neq 0$ for any $I\in \Fac\SC(Q,\pi)$ and $i\in I$.

\begin{lemma}\label{lem:regular-functional}
Suppose that there exists a regular functional $f_0$ such that there is a unique facet $I$ such that $f_0(u_i^I)>0$ for all $i\in I$.
Then for any regular functional $f$ in $V^*$, there is a unique facet $J$ such that $f(u_i^J)>0$ for all $i\in J$.
\end{lemma}

\begin{proof}
The finite arrangement of distinct hyperplanes $(u_i^I)^\perp$ gives a chamber decomposition $\bigsqcup_{t\in T} C_t$ of the set of regular functionals on $V$. Here $T$ is the index set of chambers in this arrangement. Fix $t\in T$. For any $f\in C_t$, we set  
\[ \mathbf{Pos}(f)= \{I\in \Fac\SC(Q,\pi) \mid f(u_i^I)>0 \text{ for all }i\in I\}.\]
The cardinality of $\mathbf{Pos}(f)$ is constant inside the chamber $C_t$. We denote it by $N(C_t)$. We only need to prove that $N(C_t)$ is unchanged after crossing any hyperplane in this finite arrangement.

Let $C_-$ and $C_+$ be two adjacent chambers separated by a hyperplane $H$, and choose $f_*$ in the relative interior of their common wall lying on no other hyperplane of this arrangement. 

Fix a facet $I\in \Fac\SC(Q,\pi)$.

If $f_*(u_i^I)\neq 0$ for all $i\in I$, there exists a neighbourhood $\mathcal N_I(f_*)$ of $f_*$ such that either $I\in \mathbf{Pos}(f)$ for all $f\in \mathcal N_I(f_*)\cap(C_-\cup C_+)$ or $I\notin \mathbf{Pos}(f)$ for all $f\in \mathcal N_I(f_*)\cap(C_-\cup C_+)$. 

If $f_*(u_i^I)=0$ for some $i\in I$, then such an $i$ is unique by the construction of $f_*$, since
$(u_i^I)_{i\in I}$ is a basis of $V$. Let $J$ be the unique facet of $\SC(Q,\pi)$ obtained by flipping $i$ to $j$. By Equation~\ref{eq:abstract-transport}, we have
\[
f_*(u_k^I)=f_*(u_k^J)\neq 0 \text{ for all } k\in I\cap J.
\]
Since $u^I_i=-u^J_j$ we have $f_*(u_j^J)=0$. If 
\[
f_*(u_k^I)=f_*(u_k^J)> 0 \text{ for all } k\in I\cap J,
\]
there exists a neighbourhood $\mathcal N_{I,J}$ of $f_*$ such that one of $I$ and $J$ is contained in $\mathbf{Pos}(f)$ for all $f\in \mathcal N_{I,J}(f_*)\cap C_-$ and the other one is contained in $\mathbf{Pos}(f)$ for all $f\in \mathcal N_{I,J}(f_*)\cap C_+$. Otherwise, there exists a neighbourhood $\mathcal N_{I,J}$ of $f_*$ such that $I,J\notin \mathbf{Pos}(f)$ for all $f\in \mathcal N_{I,J}(f_*)\cap(C_-\cup C_+)$.

Since $\SC(Q,\pi)$ has finitely many facets, there are only finitely many neighbourhoods $\mathcal N_I(f_*)$ and $\mathcal N_{I,J}(f_*)$
as above. Let $\mathcal N$ be their intersection. Then $\mathcal N$ is an open neighbourhood of $f_*$.

Since $f_*$ lies in the closure of both $C_-$ and $C_+$, we can choose
\[
f_-\in \mathcal N\cap C_- \text{ and } f_+\in \mathcal N\cap C_+.
\] We have
\[
|\mathbf{Pos}(f_-)|= |\mathbf{Pos}(f_+)|.
\]
Hence we have $N(C_-)=N(C_+)$.

Since any two chambers can be joined by a sequence of adjacent pairs,
the numbers $N(C_t)$ are equal for all $t\in T$. Since $|\mathbf{Pos}(f_0)|=1$, we have $N(C_t)=1$ for all $t\in T$.
\end{proof}

\begin{proposition}\label{prop:convexity-criterion}
Assume the hypotheses of Lemma~\ref{lem:regular-functional}.
Suppose that there exists a family of points $b_I\in V$ indexed by facets of $\SC(Q,\pi)$ satisfying
\begin{equation}\label{eq:abstract-edges}
b_J-b_I=\mu_{I,i}u_i^I,\qquad \mu_{I,i}>0
\end{equation}
at every flip $I\setminus\{i\}=J\setminus\{j\}$. Then $P=\conv\{b_I:I\in\Fac\SC (Q,\pi)\}$ is a simple polytope with tangent cones
\begin{equation}\label{eq:tangent-cone}
    \cone_{b_I}P=\cone\{u_i^I:i\in I\}.
\end{equation}

Its nonempty faces are indexed bijectively by faces of $\SC(Q,\pi)$:
\[
P_F=\conv\{b_I:F\subseteq I\in\Fac\SC (Q,\pi)\},\qquad F\in\SC (Q,\pi).
\]
\end{proposition}

\begin{proof}
We first prove that the vertex set of $P$ is just the set $B=\{b_I\}_{I\in \Fac\SC(Q,\pi)}$. Since $P=\conv B$, the vertex set of $P$ is a subset of $B$. For any regular functional $f$, there exists a facet $I$ such that $f(b_I)=\min\{f(x)\mid x\in P\}$. By Lemma~\ref{lem:regular-functional}, such a facet is unique. Thus $\mu_{I,i}f(u^I_i)=f(b_J)-f(b_I)\geq 0$ for any flip $I \overset{\{i,j\}}{\longrightarrow} J$ and any $i\in I$. Since $\mu_{I,i}>0$ and $f$ is regular, we have $f(u^I_i)>0$ for all $i\in I$. By Lemma~\ref{lem:regular-functional}, different facets of $\SC(Q,\pi)$ give different points in $B$.

For any facet $I$, since $\{u^I_i\}_{i\in I}$ is a basis of $V$, we can choose a regular functional $f$ such that $f(u^I_i)>0$ for any $i\in I$. Precisely, we first define a functional $g$ by $g(u^I_i)=1$ for all $i\in I$ and then perturb it to get such a regular functional $f$. The argument above implies that there exists a unique facet $K$ such that $f(b_K)=\min\{f(x)\mid x\in P\}$ and $f(u^K_i)>0$ for all $i\in K$. By Lemma~\ref{lem:regular-functional}, we have $K=I$. Thus, $b_I$ is a vertex of $P$. Therefore, the set $B$ is just the vertex set of $P$. 

For every $I$, set
\[
C=\cone_{b_I}P=\{\sum_{J\in \Fac\SC(Q,\pi)}t_J(b_J-b_I)\mid t_J\geq0\},
\]
and 
\[
D=\cone\{u_i^I:i\in I\}=\{\sum_{i\in I}\lambda_iu^I_i\mid \lambda_i\geq 0\}.
\]
Equation~\eqref{eq:abstract-edges} gives $D\subseteq C$.

Suppose $y\in C\setminus D$. Since $\{u^I_i\}_{i\in I}$ is a basis of $V$, $y=\sum_i a_iu_i^I$ and there exists $i_0\in I$ such that $a_{i_0}<0$. Thus there exists a regular functional $f$ such that $f(u_i^I)>0$ for any $i\in I$ and $f(u^I_{i_0})>-\frac{1}{a_{i_0}}(\sum_{i\neq i_0, i\in I}a_if(u^I_i))$. We have $f(y)<0$. Since $f(u_i^I)>0$ for any $i\in I$, we know that $f(b_J)>f(b_I)$ for all facets $J\neq I$. We have $f(y)\geq 0$. Contradiction. Therefore $C=D$.

Since $B$ is the vertex set of $P$, for any vertex $b_I$ of $P$, the $d$ independent rays of $D$ given by flips are precisely all $d$ edges adjacent to $b_I$. Thus $P$ is a simple $d$-dimensional polytope.

For any face $F$ of $\SC (Q,\pi)$, choose a facet $I\supseteq F$ and define a functional $f_F$ by
\[
f_F(u_i^I)=
\begin{cases}1,&i\in F,\\0,&i\in I\setminus F.\end{cases}
\]
Equation~\ref{eq:tangent-cone} implies 
\[
P_F=\{x\in P:f_F(x-b_I)=0\}
\]
is a face of $P$. Since the facets containing $F$
are connected by flips retaining $F$, we have $P_F\supseteq\conv\{b_I:F\subseteq I\in\Fac\SC (Q,\pi)\}$.

Conversely, let $b_J$ be a vertex of $P_F$. Since the edge
graph of $P_F$ is connected, choose an edge path
\[
b_I=b_{K_0},b_{K_1},\ldots,b_{K_r}=b_J
\]
inside $P_F$. We prove inductively that
\[
F\subseteq K_t,\text{ and } f_F(u_i^{K_t})=1, i\in F \text{ and } f_F(u_i^{K_t})=0, i\in K_t\setminus F.
\]
This holds for $K_0=I$ by definition. Suppose that it holds for
$K_t$ and $K_t\overset{\{i,j\}}{\longrightarrow}K_{t+1}$, where $i\in K_t$.
Since $b_{K_t}$ and $b_{K_{t+1}}$ lie in $P_F$,
\[
0=f_F(b_{K_{t+1}}-b_{K_t})=\mu_{K_t,i}f_F(u_i^{K_t}).
\]
Thus $f_F(u_i^{K_t})=0$, so $i\notin F$ and $F\subseteq K_{t+1}$.
For every $l\in K_{t+1}\setminus \{j\}$, Equation~\ref{eq:abstract-transport} gives
\[
f_F(u_l^{K_{t+1}})
=f_F(u_l^{K_t}).
\]
Moreover, $j\notin K_t$ and $F\subseteq K_t$, so $j\notin F$. We have
\[
f_F(u_j^{K_{t+1}})
=-f_F(u_i^{K_t})=0.
\]
This completes the induction and shows that $F\subseteq J$. Thus we have
\[
\Vertices(P_F)
=\{b_J:F\subseteq J\in\Fac\SC(Q,\pi)\}.
\]

Finally, let $G$ be a nonempty face containing a vertex $b_I$. Choose $g\in V^*$ such that $G=\{x\in P:g(x-b_I)=0\}$ and $g(x-b_I)\geq0 $ for all $x\in P$. Set $F=\{i\in I:g(u_i^I)>0\}$. We have $g(u_i^I)\geq0$ for all $i\in I$, and every $x\in P$ satisfies
$x-b_I=\sum_{i\in I}a_i u_i^I$ with $a_i\geq0$. Thus we have
\[
g(x-b_I)=0 \iff a_i=0\ (i\in F) \iff f_F(x-b_I)=0.
\]
Therefore $G=P_F$. Moreover, the tangent cone of $P_F$ at $b_I$
is generated by $\{u_i^I\}_{i\in I\setminus F}$.
Since $\{u_i^I\}_{i\in I}$ is a basis, such a subset $F$ of $I$ is unique. Hence the nonempty faces of $P$ are indexed bijectively by the faces of $\SC(Q,\pi)$.
\end{proof}

\begin{corollary}\label{cor:polytopal}
Under the hypotheses of Proposition~\ref{prop:convexity-criterion},
the subword complex $\SC(Q,\pi)$ is polytopal. For any $z\in\operatorname{int}(P)$, the boundary complex of the polar dual polytope $(P-z)^\circ$ is isomorphic to $\SC(Q,\pi)$ as simplicial complex.
\end{corollary}

\begin{proof}
The inclusion-reversing face correspondence $F\mapsto P_F$ gives the required isomorphism.
\end{proof}

\begin{remark}
    The cone equality and the face argument above are similar to the proof given in \cite[Theorem~4.4 and Proposition~4.6]{JS}. Homever, we work in an arbitrary vector space $V$, which need not be the root space. 
\end{remark}

\begin{proof}[Proof of Theorem~\ref{thm:km-ri}]
Apply the convexity criterion in the space $E_Q$, with $u_i^I=-r(I,i)$ and $b_I=B(I)-B(I^+)$.
The flip graph is connected and all brick-vector
differences along flips lie in $E_Q$, so $b_I\in E_Q$.
Notice that $(u_i^I)_{i\in I}$ is a basis of $E_Q$ for every facet $I$ since $\SCQ$ is root independent.

Choose $h\in V^*$ by $h(\alpha_s)=1$ for all $s\in S$,
where $V$ is the root space, and set $f_0=-h|_{E_Q}$.
Since any root is either positive or negative, the functional $f_0$ is regular on the vectors $(u_i^I)_{i\in I}$. We have $f_0(u_i^I)>0\text{ for all }i\in I$ if and only if $R(I)\subseteq\Phi^+$. The unique such facet is $I^+$. Proposition~\ref{prop:convexity-criterion} and Corollary~\ref{cor:polytopal} give the result.
\end{proof}

The proof of Theorem~\ref{thm:km-ri} can be generalized to an arbitrary Coxeter system. We use the real reflection realization in Section~\ref{sec:coxeter}. Here the coefficient $a_{I,i,k}$ in Lemma~\ref{lem:weight-flip} is a real number instead of an integer. Proposition~\ref{prop:convexity-criterion} can be applied directly.
\begin{proposition}\label{prop:general-realization}
    For an arbitrary Coxeter system $(W,S)$, any spherical root independent subword complex is polytopal.
\end{proposition}

\section{Double root free case}\label{sec:double}
Now we deal with the double root free case. We assume that $\SC (Q,\pi)$ is double root free throughout this section. We first deal with the Kac-Moody case using the lattice language and then explain that the same argument can be applied to an arbitrary Coxeter system by using the vector space instead of the integral lattice. 

The idea is to create a lattice $M_Q$ such that the root configuration space $E_Q$ is a quotient of $M_Q\otimes\mathbb{R}$ in order to construct a basis in $M_Q\otimes\mathbb{R}$ for each facet $I$ and apply Proposition~\ref{prop:convexity-criterion}.

For each facet $I$, let $M_I=\bigoplus_{i\in I}\Z e_i^I$, and construct a map $p_I$ from $M_I$ to $V$ by $p_I(e_i^I)=-r(I,i)$.

Consider a flip $I\setminus\{i\}=J\setminus\{j\}$ and set $\beta=r(I,i)$. For $k\in I\cap J$, define
\[
c_k=
\begin{cases}
-\langle r(I,k),\beta^\vee\rangle,
&\min(i,j)<k<\max(i,j),\\
0,&\text{otherwise}.
\end{cases}
\]
In particular, $c_k\in\Z$. Lemma~\ref{lem:flip} gives
\[
r(J,j)=-r(I,i) \text{ and } r(J,k)=r(I,k)+c_k r(I,i) \text{ for all } k\in I\cap J.
\]
We identify $M_I$ and $M_J$
by
\begin{equation}\label{eq:fold-transport}
e_j^J=-e_i^I \text{ and } e_k^J=e_k^I+c_k e_i^I \text{ for all } k\in I\cap J.
\end{equation}
These formulas define an isomorphism of lattices and the maps $p_I$ and $p_J$ are compatible with this identification.

\begin{proposition}\label{prop:fold-lattice}
The identifications \eqref{eq:fold-transport} define a lattice $M_Q\simeq M_I$ of rank $d$, with
a basis $(e_i^I)_{i\in I}$ for every facet $I$.
The maps $p_I$ define an exact sequence
\begin{equation}\label{eq:fold-exact}
0\longrightarrow K_Q^{\Z} \longrightarrow M_Q \overset{p}{\longrightarrow}\Lambda_Q\longrightarrow0,
\end{equation}
where $\Lambda_Q=\sum_{i\in I}\Z r(I,i)$ is a lattice in $V$ independent of the facet $I$.
\end{proposition}

\begin{proof}
Without loss of generality, we assume that the cardinality $d$ of any facet of $\SC(Q,\pi)$ is larger than 2. In fact, if $d\leq 2$, any double root free subword complex is root independent by definition. The flip graph is the
one-skeleton of the cell decomposition dual to the subword complex. A spherical subword complex is a sphere. Since the sphere has dimension $d-1\ge2$, it is simply connected. Thus every closed path in the flip graph is generated by cycles coming from one-dimensional links in the subword complex. Thus we only need to check that the identifications are compatible for cycles induced by one-dimensional links.

Let $F$ be a face with $|F|=d-2$. Its link $L$ is a cycle of facets. Choose facet $I=F\sqcup\{i,j\}$ on this cycle and set $U_I=\Span_{\R}\{e_i^I,e_j^I\}$. Along the cycle, every flip retains $F$. By \eqref{eq:fold-transport}, for any facet $I'=F\sqcup\{i',j'\}$ on this cycle, the space $U_{I'}=\Span_{\R}\{e_{i'}^{I'},e_{j'}^{I'}\}$ can be identified with $U_I$. We simply denote this two dimensional space by $U$.

Let $H$ be the linear map from $M_I\otimes\mathbb{R}$ to itself obtained by going once around the cycle starting from and ending at $I$. We denote $M_I\otimes\mathbb{R}$ by $M_{I,\mathbb{R}}$ and identify $p_I$ with its induced linear map from $M_{I,\mathbb{R}}$ to $V$.  We have $(H-\mathrm{id})M_{I,\R}\subseteq U$ and $p_IH=p_I$.
By double root freeness, $r(I,i)$ and $r(I,j)$ are not proportional in $\Lambda_Q$. Hence $p_I$ is injective on $U$. The inclusion and the equality above imply $H=\mathrm{id}$.

We can therefore identify all $M_I$ with a single lattice $M_Q$. The maps $p_I$ agree with these identifications and hence define a map $p$ on $M_Q$. Equation~\eqref{eq:fold-transport} gives an integral change of bases between the bases $\{e^I_i\}_{i\in I}$ for $I\in \Fac\SC(Q,\pi)$. This common image is $\Lambda_Q$, giving \eqref{eq:fold-exact}.
The kernel $K_Q^{\Z}$ is free, being a subgroup of the free abelian group $M_Q$. 
\end{proof}

\begin{example}
In type $A_1$, let $Q=sss$, $\pi=s$. The subword complex is the boundary of a triangle, which is not double root free.
Start with the free abelian group $M_{\{1,2\}}=\Z a\oplus\Z b$ and
$p(a)=p(b)=-\alpha$. We consider a cycle of flips:
\[
\{1,2\}\overset{\{2,3\}}{\longrightarrow}\{1,3\}\overset{\{1,2\}}{\longrightarrow}\{2,3\}\overset{\{3,1\}}{\longrightarrow}\{1,2\}.
\]
The map induced by the flip rule around this cycle is not the identity. Precisely, it has the matrix form
\[
H=\begin{pmatrix}0&-1\\1&2\end{pmatrix}\neq I_2.
\]
Hence $M_Q$ is not well defined in this case. Thus we can not identify the lattice $M_I$ for all facets $I$ in this case.
\end{example}
The root configuration space $E_Q$ is just $\Lambda_Q\otimes\mathbb{R}$. We have a surjective map $p$ from $M_{Q,\R}$ to $E_Q$ induced by the maps $p_I$. Now we lift the weight functions of $\SC(Q,\pi)$ to vectors in $M_{Q,\R}$:

\begin{proposition}\label{prop:lifted-weights}
There are unique \textbf{lifted brick weights} $\wt w(I,k)\in M_{Q,\R}$ satisfying
\begin{equation}\label{eq:lifted-weight-rule}
\wt w(I^+,k)=0 \text{ and } \wt w(J,k)-\wt w(I,k)=a_{I,i,k}e_i^I
\end{equation}
for any flip $I\setminus\{i\}=J\setminus\{j\}$ and $k\in [m]$. They satisfy
\[
p(\wt w(I,k))=w(I,k)-w(I^+,k).
\]
\end{proposition}

\begin{proof}
We first normalize $\wt w(I^+,k)=0$ for all $k\in [m]$. Fix $k$. For each facet $I$, choose a flip path
\[
I^+=I_0,I_1,\ldots,I_N=I,
\]
where we have flips $I_t\setminus\{i_t\}=I_{t+1}\setminus\{j_{t+1}\}$ and $i_t\in I_t$. We define
\[
\wt w(I,k)
=\sum_{t=0}^{N-1}a_{I_t,i_t,k}e_{i_t}^{I_t}\in M_Q.
\]
We need to check that this sum does not depend on the chosen path.

As in the argument in Proposition~\ref{prop:fold-lattice}, it suffices to check that the equation 
\[
z=\sum_t\bigl(\wt w(I_{t+1},k)-\wt w(I_t,k)\bigr)=0.
\]
holds whenever the flip path is a cycle for a one-dimensional link with common face $F$. Denote the two-dimensional subspace of $M_{Q,\R}$ defined by this face $F$ as in the proof of Proposition~\ref{prop:fold-lattice} by $U$. We have $z\in U$ and
\[
p(z)=-\sum_t a_{I_t,i_t,k}r(I_t,i_t)=\sum_t\bigl(w(I_{t+1},k)-w(I_t,k)\bigr)=0.
\]
Since $p|_{U}$ is injective by double root freeness, we obtain $z=0$. Hence the definition is independent of the path.

Finally, the normalization condition for $\wt w(I^+,k)$, $k\in[m]$, and the flip rule give the uniqueness part.
\end{proof}

Now we can define \textbf{lifted brick vectors} and the $\textbf{lifted brick polytope}$ in $M_{Q,\R}$ as follows:
\begin{equation}\label{eq:lifted-brick}
\wt B(I)=\sum_{k=1}^m\wt w(I,k) \text{ and } P_Q=\conv\{\wt B(I):I\in\Fac\SC (Q,\pi)\}.
\end{equation}
By Lemma~\ref{lem:weight-flip} and Proposition~\ref{prop:lifted-weights}, we have
\begin{equation}\label{eq:lifted-brick-flip}
\wt B(J)-\wt B(I)=\lambda_{I,i}e_i^I \text{ where } \lambda_{I,i}>0.
\end{equation}

\begin{theorem}\label{thm:drf-polytope}
For a double root free spherical subword complex of Kac-Moody type, $P_Q$ is a simple $d$-dimensional lattice
polytope. Its vertices are the distinct lifted brick vectors $\wt B(I)$, its edges are given by flips and
\begin{equation}\label{eq:lifted-tangent}
\cone_{\wt B(I)}P_Q=\cone\{e_i^I:i\in I\}.
\end{equation}
Its nonempty faces are
\[
P_F=\conv\{\wt B(I):F\subseteq I\in\Fac\SC (Q,\pi)\}.
\]
The boundary of the polar of an interior translate of $P_Q$ gives the polytopal realization of $\SC (Q,\pi)$. If $\SC (Q,\pi)$ is root independent, $p_{\R}:M_{Q,\R}\to E_Q$ is an isomorphism and this agrees with Theorem~\ref{thm:km-ri}.
\end{theorem}
\begin{proof}
Use Proposition~\ref{prop:convexity-criterion} with the space $M_{Q,\R}$, $u_i^I=e_i^I$, and $b_I=\wt B(I)$. Equations~\eqref{eq:fold-transport} and Equation~\eqref{eq:lifted-brick-flip} induce the conditions \eqref{eq:abstract-transport} and \eqref{eq:abstract-edges}. We take a functional $h\in V_{}^*$ such that $h|_{\Phi^+}>0$. Take $f=-h\circ p$. This is a regular functional in $M_{Q,\R}^*$ such that $f(e^{I^+}_i)>0$ for all $i\in I^+$. Lemma~\ref{lem:regular-functional} therefore holds.

Proposition~\ref{prop:convexity-criterion} and Corollary~\ref{cor:polytopal} give the desired assertion. 
\end{proof}

\begin{remark}\label{rem:general-Coxeter}
In fact, the polytopal realization of the spherical double root free subword complex can be generalized to an arbitrary Coxeter group using the real reflection realization in Section~\ref{sec:coxeter} and apply Proposition~\ref{prop:convexity-criterion}. We construct a lifted space of the root configuration space $E_Q$ directly and obtain an exact sequence of vector space which is which is an analogue of Proposition~\ref{prop:fold-lattice}. Then apply Proposition~\ref{prop:convexity-criterion} to obtain Theorem~\ref{thm:drf-polytope}.

The weight functions, the brick vectors and the lifted brick vectors can be defined in the same way. But in general we can only get a polytope in a real vector space rather than a lattice. 
\end{remark}

\section{Relations with braid varieties}
\label{sec:complex}
We use braid variety and the set of unfolded galleries in buildings to interpret $M_Q$ in the Kac-Moody case. First we fix some building conventions following \cite{AB}.

The Coxeter complex $\Sigma(W,S)$ is the chamber complex whose chambers $C_w$ are indexed by $w\in W$, with $C_w$ and $C_{ws}$ sharing a panel of type $s$. A building $X$ of type $(W,S)$ is a chamber complex covered by apartments isomorphic to $\Sigma(W,S)$. Any two chambers lie in an apartment, and any two apartments are related by a type-preserving isomorphism fixing their intersection. It is thick if each panel belongs to at least three chambers. The Weyl distance $\delta(D,E)\in W$ is computed in an apartment containing $D,E$; in the standard apartment, we have
\[
 \delta(C_u,C_v)=u^{-1}v.
\]

A folded gallery of type $Q=q_1\cdots q_m$ is a sequence
$(D_0,\ldots,D_m)$ with
$\delta(D_{k-1},D_k)\in\{e,q_k\}$.
The $k$-th step is a fold if $D_{k-1}=D_k$ and a crossing otherwise. The gallery is unfolded if every step is a crossing; an unfolded gallery is minimal if $m=\ell(\delta(D_0,D_m))$.
In the standard apartment an unfolded gallery
of type $Q$ starting at $C_e$ is unique and ends at
$C_{q_1\cdots q_m}$.

Given a word $Q$, a set $J\subseteq[m]$ determines a gallery
$\gamma_J=(C_{x_0^J},\ldots,C_{x_m^J})$ of type $Q$ with folded positions $J$. It ends at $C_\pi$ exactly when $\sigma_J=\pi$. 

For the algebraic constructions below, assume that $\SC(Q,\pi)$
is a spherical double root free subword complex. Set $d=m-\ell(\pi)$.
Let $G$ be a semisimple simply connected complex group in finite Weyl type, or a minimal split Kac-Moody group over $\C$ in Kac-Moody type. Write $H$ for its Cartan torus, $B$ for its positive Borel subgroup, and $P_s$ for its minimal standard parabolic subgroups. For each real root $\alpha$, let $U_\alpha$ be the corresponding root subgroup of $G$. We take $X$ to be the building with chambers $G/B$, and standard apartment $A$ with
chambers $C_w=\dot wB$. In Kac-Moody type this is the positive half of the twin building, which is a building itself \cite{AB}. Its panel residues are $gP_s/B\simeq\mathbb P^1$, and
\[
 \delta(gB,hB)=w\quad\Longleftrightarrow\quad
 g^{-1}h\in B\dot wB.
\]
For each real root $\alpha$, choose
$x_\alpha:(\C,+)\xrightarrow{\sim}U_\alpha$. We use the standard
identities
\begin{equation}\label{eq:root-group-actions}
 \dot wU_\alpha\dot w^{-1}=U_{w\alpha},\qquad
 h x_\alpha(t)h^{-1}=x_\alpha\bigl(\alpha(h)t\bigr)
 \text{ for all } h\in H.
\end{equation}
Here roots are regarded as characters of $H$. The pairing of
characters and cocharacters is defined by
$\lambda(\nu(z))=z^{\langle\lambda,\nu\rangle}$ and agrees with the pairing in Section~\ref{sec:coxeter}. See \cite{Kumar} for more details.

The Bott-Samelson variety
\[
 \operatorname{BS}(Q)=P_{q_1}\times^B\cdots\times^B P_{q_m}/B
\]
parametrizes all galleries of type $Q$ starting at $C_e$. Let $Z_Q$ be the fiber of its endpoint map over $C_\pi$ in \cite{BGS}.

We recall standard facts about panels and retractions in
\cite{AB}. For $s\in S$, let $P_s=B\sqcup B\dot{s}B$ be the
corresponding standard minimal parabolic subgroup. A panel residue $\mathcal R=gP_s/B$ has a unique gate $G_{\mathcal R}$ towards a chamber $D$, satisfying
\[
 \operatorname{dist}(D,E)
 =\operatorname{dist}(D,G_{\mathcal R})
   +\operatorname{dist}(G_{\mathcal R},E)
 \qquad(E\in\mathcal R),
\]
where $\operatorname{dist}(D,E)=\ell(\delta(D,E))$.
For $z\in W$, the retraction $\rho_z:X\to A$ centered at $C_z$ fixes $A$ and restricts to an isomorphism on
each apartment containing $C_z$. It sends the gate towards $C_z$ of a panel residue to one of the two chambers of its image panel, and all other chambers of the residue to the other.  If a type-$s$ step in a unfolded gallery retracts to a fold at $C_w$ by $\rho_z$, we have $z^{-1}w\alpha_s\in\Phi^-$.

In the split complex Kac-Moody case, if $\beta=w\alpha_s$,
the root group $U_{-\beta}$ fixes $C_{ws}$ and acts simply
transitively on the complement of $C_{ws}$ in their common panel residue. We have a parameterization
\[
 \C\xrightarrow{\sim}\mathcal R\setminus\{C_{ws}\},
 \qquad t\longmapsto x_{-\beta}(t)C_w,
\]
with $t=0$ giving $C_w$. We keep $M_Q$ as the integral lattice of Proposition~\ref{prop:fold-lattice}. Recall that $I^+$ is the unique facet of $\SCQ$ such that the root configuration contains only positive roots.

Now, we give an algebraic interpretation of the lattice $M_Q$ in the complex case.
\begin{definition}\label{def: com-braid-variety}
    Fix an arbitrary Coxeter system $(W,S)$. Given a word $Q$ whose Demazure product is $\pi\in W$ and a thick building $X$ of type $(W,S)$, \textbf{the set of unfolded galleries} of type $Q$ from $C_e$ to $C_{\pi}$ in $X$ is denoted by
    \begin{align}
        \mathcal{G}(Q):=\{ \gamma: C_e \longrightarrow C_{\pi} \mid \gamma \text{ is an unfolded gallery of type } Q \text{ in } X\}.
    \end{align}
\end{definition}
In the finite case, $\mathcal G(Q)$ is just the braid variety as considered in \cite{BGS}.

\begin{proposition}\label{prop:complex-torus}
Assume that $\SC(Q,\pi)$ is double root free and $X=G/B$. Then we have a parameterization of $\mathcal{G}(Q)$ given by $I^+$ of $\SCQ$:
\[
\mathcal G(Q)\simeq(\C^\times)^{|I^+|}.
\]
Writing $t_i$ for the coordinate at $i\in I^+$, the assignment $e_i^{I^+}\mapsto[t_i]$ gives an isomorphism from $M_Q$ to the character lattice of $\mathcal G(Q)$:
\begin{equation}\label{eq:character-identification}
M_Q\xrightarrow{\sim}\mathcal O(\mathcal G(Q))^\times/\C^\times.
\end{equation}
Set $T_Q=\Spec\C[M_Q]$. The Cartan action is translation through
$\eta:H\to T_Q$, with character map
\[
\eta^*:M_Q\xrightarrow{\,p\,}\Lambda_Q\hookrightarrow X^*(H),
\qquad e_i^{I^+}\longmapsto-r(I^+,i).
\]
\end{proposition}

\begin{proof}
We set $v_k=x_k^{I^+}$ and $\beta_k=r(I^+,k)=v_{k-1}\alpha_{q_k}$.
As in the proof of Lemma~\ref{lem:chamber-maps-for-e}, we have $R(I^+)\subseteq \Phi^+\cap\pi\Phi^-$. Thus, we have $\beta_k\in \Phi^+\cap\pi\Phi^-$ for all $k\in[m]$.

Let $\mathcal R_k$ be the type-$q_k$ panel residue containing $C_{v_{k-1}}$. Its gate towards $C_\pi$ is therefore $G_k=C_{v_{k-1}q_k}$. We have that $U_{-\beta_k}$ fixes $C_\pi$ since $\dot\pi^{-1}U_{-\beta_k}\dot\pi
 =U_{-\pi^{-1}\beta_k}\subset B$.
 
The root group action gives an algebraic isomorphism
\[
 \phi_k:\mathbb C\xrightarrow{\sim}
 \mathcal R_k\setminus\{G_k\},
 \qquad t\longmapsto x_{-\beta_k}(t)C_{v_{k-1}},
\]
where $\phi_k(0)=C_{v_{k-1}}$. Its restriction to $\C^\times$ parametrizes $\mathcal R_k\setminus\{C_{v_{k-1}},G_k\}$.

Now we give a parametrization of $\mathcal G(Q)$ using the retraction $\rho_\pi$ of $X$ onto the standard apartment centered at $C_\pi$.

We first show that $\rho_{\pi}(\mathcal G(Q))=\{\gamma_{I^+}\}$.
Suppose that $\rho_\pi(\gamma)=\gamma_J$. The endpoint of $\gamma_J$ is $C_\pi$, and so the product of the complementary subword of $J$ is $\pi$. By Theorem~\ref{thm:double}, this subword is a reduced expression of $\pi$ and $J$ is
a facet of $\SCQ$. At a folded position $j\in J$, the two distinct chambers have the same image under $\rho_\pi$, so neither is the gate of their common panel residue. Thus, we have
\[
 \pi^{-1}r(J,j)\in\Phi^-.
\]
However, the condition $|r(J,j)|\in\Phi^+\cap\pi\Phi^-$ implies that $r(J,j)\in\Phi^+$. We have $R(J)\subseteq\Phi^+$. By the uniqueness in Lemma~\ref{lem:chamber-maps-for-e}, we have $J=I^+$.

Given $t=(t_i)_{i\in I^+}\in(\C^\times)^{I^+}$, we set $g_0=1$,
$D_0=C_e$, and define inductively
\[
 g_k=
 \begin{cases}
 g_{k-1},&k\notin I^+,\\
 g_{k-1}x_{-\beta_k}(t_k),&k\in I^+,
 \end{cases}
 \qquad D_k=g_kC_{v_k}.
\]
We verify that these chambers form an unfolded gallery.
Suppose that $g_{k-1}$ fixes $C_\pi$ and
$D_{k-1}=g_{k-1}C_{v_{k-1}}$. The type-$q_k$ panel residue
containing $D_{k-1}$ is $g_{k-1}\mathcal R_k$, whose gate
towards $C_\pi$ is $g_{k-1}G_k$.

If $k\notin I^+$, we have $v_k=v_{k-1}q_k$ and
\[
 D_k=g_{k-1}G_k\neq D_{k-1}.
\]

If $k\in I^+$, we have $v_k=v_{k-1}$ and
\[
 D_k=g_{k-1}\phi_k(t_k)
 \in g_{k-1}\mathcal R_k
 \setminus\{D_{k-1},g_{k-1}G_k\}.
\]
In both cases, $D_{k-1}$ and $D_k$ are distinct and
$q_k$-adjacent chambers in $X$. Also, $g_k$ fixes $C_\pi$, because every root group $U_{-\beta_k}$ fixes $C_\pi$. This completes the induction.

Since $v_m=\pi$ and $g_m$ fixes $C_\pi$, we have
$D_m=g_mC_\pi=C_\pi$.
Each $g_k(t)$ is a fixed ordered product of root group
elements $x_{-\beta_i}(t_i)$ with $i\in I^+$ and $i\leq k$,
and therefore depends regularly on $t$.
Since the root group actions are regular, each
$t\mapsto D_k(t)=g_k(t)C_{v_k}$ is regular.
Thus, the construction defines a morphism
\[
\Theta:(\C^\times)^{I^+}\longrightarrow\mathcal G(Q).
\]
Now we consider the inverse direction to prove that $\Theta$ is an isomorphism. Let $\gamma=(D_0,\ldots,D_m)\in\mathcal G(Q)$.
We construct its parameters inductively. Suppose that
$g_{k-1}$ has been defined, fixes $C_\pi$, and satisfies
$D_{k-1}=g_{k-1}C_{v_{k-1}}$.
We have $\rho_\pi(\gamma)=\gamma_{I^+}$.

If $k\notin I^+$, we have
\[
 D_k=g_{k-1}G_k.
\]
We then set $g_k=g_{k-1}$.

If $k\in I^+$, we have
\[
 g_{k-1}^{-1}D_k
 \in\mathcal R_k\setminus\{C_{v_{k-1}},G_k\}.
\]
There is therefore a unique $t_k=\phi_k^{-1}(g_{k-1}^{-1}D_k)\in\C^\times$ such that $g_k=g_{k-1}x_{-\beta_k}(t_k)$.

Inductively, we can reconstruct all parameters of $\gamma$. We have a map
\[
 \Psi:\mathcal G(Q)\longrightarrow(\C^\times)^{I^+}.
\]
To prove that $\Psi$ is regular, we proceed inductively
over the positions in $I^+$. Suppose that the parameters $t_i$, for $i\in I^+$ with $i<k$, are already regular functions on $\mathcal G(Q)$.
We have that $g_{k-1}(\gamma)$ is a fixed ordered product of
root group elements with regular parameters.

Notice that each evaluation map
\[
\operatorname{ev}_k:\operatorname{BS}(Q)\longrightarrow G/B,
\qquad
(D_0,\ldots,D_m)\longmapsto D_k,
\]
is the restriction of a coordinate projection and hence is
an algebraic morphism. Since $\mathcal G(Q)$ carries the open subvariety structure inside $\operatorname{ev}_m^{-1}(C_\pi)$, the restriction of $\operatorname{ev}_k$ to $\mathcal G(Q)$ is also a morphism.

The evaluation map $\gamma\mapsto D_k$ and the root group actions are regular. The inverse is given by $x_\alpha(a)^{-1}=x_\alpha(-a)$.
Thus, the map
\[
F_k:\mathcal G(Q)\longrightarrow
\mathcal R_k\setminus\{C_{v_{k-1}},G_k\},
\qquad
\gamma\longmapsto g_{k-1}(\gamma)^{-1}D_k
\]
is regular. Since $\phi_k^{-1}$ is regular on
$\mathcal R_k\setminus\{G_k\}$, the function $t_k=\phi_k^{-1}\circ F_k$ is regular and takes values in $\C^\times$. Starting with $g_0=1$, this proves inductively that all coordinates of $\Psi$ are regular, hence $\Psi$ is a morphism. The two constructions are inverse to each other, so $\Theta$ is an algebraic isomorphism.

Now we can write
\[
 \mathcal O(\mathcal G(Q))
 =\C[t_i^{\pm1}:i\in I^+],
 \qquad
 \mathcal O(\mathcal G(Q))^\times
 =\left\{c\prod_{i\in I^+}t_i^{n_i}:
 c\in\C^\times,\ n_i\in\Z\right\}.
\]
Since $(e_i^{I^+})_{i\in I^+}$ is a basis of $M_Q$, the map
\[
 \sum_{i\in I^+}n_i e_i^{I^+}
 \longmapsto
 \left[\prod_{i\in I^+}t_i^{n_i}\right]
\]
gives Equation~\eqref{eq:character-identification}.

Finally, $H$ acts on $\mathcal G(Q)$ by $h\cdot(D_0,\ldots,D_m)=(hD_0,\ldots,hD_m)$.

For any chamber $C_w=\dot wB$ in the standard apartment, we have $ h\dot wB=\dot w(\dot w^{-1}h\dot w)B=\dot wB$. For $h\in H$, set $t_i'=(-\beta_i)(h)t_i$. The conjugation of the root group gives $h x_{-\beta_i}(t_i)h^{-1}
 =x_{-\beta_i}(t_i')$. The induction above gives $g_k(t')=h g_k(t)h^{-1}$. Since $hC_{v_k}=C_{v_k}$, we have
\[
 g_k(t')C_{v_k}=h g_k(t)h^{-1}C_{v_k}=hg_k(t)C_{v_k}.
\]
The Cartan action can be expressed as
\[
 h\cdot(t_i)_{i\in I^+}
 =\bigl((-\beta_i)(h)t_i\bigr)_{i\in I^+}.
\]
Identifying $\mathcal G(Q)$ with $T_Q$ by $\chi^{e_i^{I^+}}=t_i$, we have a morphism $\eta:H\longrightarrow T_Q$ by requiring $\chi^{e_i^{I^+}}(\eta(h))=(-\beta_i)(h)$.
Hence
\[
 \eta^*(e_i^{I^+})=-\beta_i=p(e_i^{I^+}),
\]
as required.
\end{proof}

\input{affine_drf_examples}

\section*{Acknowledgements}

The author is very grateful to Stéphane Gaussent, Mikhail Gorsky, Paul Philippe and Petra Schwer for many useful discussions. The author acknowledges support from the Procope project "Buildings, galleries and beyond" and the PhD funding from Institut Camille Jordan, Université Jean Monnet.

\section*{Use of AI}
GPT 5.6 Sol and GPT 6.1 Sol were used to check grammar, search for references, provide \LaTeX{} code, and double-check statements after they had been checked by the author.

GPT 6 Astra provided the author with suggestions for improving the clarity of the proof of Lemma~\ref{lem:regular-functional}. GPT 6 Astra helped the author find a suitable geometric representation of an arbitrary Coxeter system in \cite{RicheLectures}. GPT 6.1 Sol helped the author to find the citable reference \cite{FU} for Equation~\eqref{eq:root-length}. GPT 6 Astra was used to compare the proof and strategy used in Section~\ref{sec:km} with those used for brick polytopes in finite case in previous work. GPT 6 Astra helped the author check the algebro-geometric aspects of the construction of the isomorphism in Proposition~\ref{prop:complex-torus}. GPT 6 Astra was used to compute examples and generate tables and figures in Section~\ref{sec:interleaved-examples}.

The author is responsible for all the mathematical proofs and constructions and he takes full responsibility for the writing of the paper. 

\bibliographystyle{plain}
\bibliography{biblio}

\bigskip
\noindent
{\footnotesize\scshape
Institut Camille Jordan UMR 5208,
Université Jean Monnet, CNRS, Centrale Lyon,
INSA Lyon, Université Claude Bernard Lyon 1,
20, rue Annino, 42023, Saint-Étienne, France
}

\par\smallskip
\noindent
{\footnotesize
\textit{Email address:}
\texttt{zijun.li@univ-st-etienne.fr}
}
\end{document}

%% file: affine_drf_examples.tex
\section{Examples of brick polytopes}
\label{sec:interleaved-examples}

In this section, we give several examples of the brick polytopes realizing the corresponding subword complex. We use the normalization $\wt B(I^+)=0$ for the lifted vectors in the space  $M_Q$. We first give a root independent example $\mathcal P_1$ of type $\widetilde A_2$. Then we give a double root free but not root independent example $\mathcal P_2$ of type $A_3$. We have $\mathcal P_2\simeq \mathcal P_1\times [0,1]$ combinatorially, although the Coxeter groups and the words are totally different. Finally, we give a double root free but not root independent example of type $\widetilde A_2$.

\subsection{Example 1}
\label{ex:int-affine}

We take $W=\widetilde A_2$, $Q=s_0s_1s_0s_1s_2s_0s_2s_0 \text{ and }
\pi=\Dem(Q)=s_0s_1s_0s_2s_0.$
The length of $\pi$ is $5$ and $d=8-5=3$. Take $I^+=\{1,2,5\}$. We have
\[
R(I^+) =(\alpha_0,\alpha_1,2\alpha_0+\alpha_1+\alpha_2)
       =(\alpha_0,\alpha_1,\alpha_0+\delta),
\text{ where } \delta=\alpha_0+\alpha_1+\alpha_2.
\]
This example is root independent. In particular, it is double root
free, and $\Lambda_Q=\Z\alpha_0\oplus\Z\alpha_1\oplus\Z\alpha_2$.

Set
\[
M_Q=\Z U\oplus\Z V\oplus\Z W \text{ and } (U,V,W)=(e_1^{I^+},e_2^{I^+},e_5^{I^+}).
\]
The projection from $M_Q$ to $\Lambda_Q$ is
\begin{equation}\label{eq:int-affine-p}
p(u,v,w)=-u\alpha_0-v\alpha_1-w(2\alpha_0+\alpha_1+\alpha_2).
\end{equation}
It is an isomorphism of lattices in this example.

All facets, coordinates of lifted brick vectors after the normalization and lifted bases are listed as follows:

\begin{table}[htbp]
\centering
\caption{$\SC_{\widetilde A_2}(s_0s_1s_0s_1s_2s_0s_2s_0,s_0s_1s_0s_2s_0)$}
\label{tab:int-base}
\medskip
\begin{tabular}{@{}ccc@{}}
\toprule
$K$ & $z_K=(u,v,w)=\wt B(K)$ & $(e_i^K)_{i\in K}$\\
\midrule
$\{1,2,5\}$ & $(0,0,0)$ & $(U,V,W)$\\
$\{1,2,8\}$ & $(0,0,2)$ & $(U,V,-W)$\\
$\{1,5,6\}$ & $(0,2,0)$ & $(U,V+W,-V)$\\
$\{1,6,7\}$ & $(0,3,1)$ & $(U,W,-V-W)$\\
$\{1,7,8\}$ & $(0,3,2)$ & $(U,-V,-W)$\\
$\{2,3,5\}$ & $(1,0,0)$ & $(U+V,-U,W)$\\
$\{2,3,8\}$ & $(1,0,2)$ & $(U+V,-U,-W)$\\
$\{3,4,5\}$ & $(2,1,0)$ & $(V,-U-V,W)$\\
$\{3,4,8\}$ & $(2,1,2)$ & $(V,-U-V,-W)$\\
$\{4,5,6\}$ & $(2,2,0)$ & $(-U,V+W,-V)$\\
$\{4,6,7\}$ & $(2,3,1)$ & $(-U,W,-V-W)$\\
$\{4,7,8\}$ & $(2,3,2)$ & $(-U,-V,-W)$\\
\bottomrule
\end{tabular}
\end{table}

For the facet $I^+$, every $w(I^+,k)$ equals $\omega_{q_k}$ except $w(I^+,6)=\omega_0-\alpha_0 \text{ and } w(I^+,8)=\omega_0-\alpha_0-\alpha_1$. The brick vector $B(I^+)$ is $4\omega_0+2\omega_1+2\omega_2-2\alpha_0-\alpha_1$.
All brick vectors are given explicitly by
\begin{equation}\label{eq:int-affine-B}
B(K)=B(I^+)-u\alpha_0-v\alpha_1-w(2\alpha_0+\alpha_1+\alpha_2),
\text{ where } (u,v,w)=z_K.
\end{equation}
For example, $B(\{1,2,8\})=B(I^+)-2(2\alpha_0+\alpha_1+\alpha_2)$.

The convex hull of the lifted vectors is
\begin{equation}\label{eq:int-R}
\mathcal P_1=
\{(u,v,w)\mid 0\leq u\leq2,\ 0\leq v\leq3,\ 0\leq w\leq2,\
u-v\leq1,\ v-w\leq2\}.
\end{equation}

\begin{figure}[htbp]
\centering

\begin{tikzpicture}[
    x={(2.2cm,0cm)},
    y={(0.95cm,0.55cm)},
    z={(0cm,1.8cm)},
    line join=round,
    line cap=round,
    facetfill/.style={fill=blue!6,draw=none},
    visedge/.style={draw=black,thick},
    hidedge/.style={draw=black!65,thick,dashed},
    vtx/.style={circle,fill=black,inner sep=1.2pt},
    vlab/.style={font=\scriptsize\bfseries,text=blue!65!black,fill=white,inner sep=1.1pt},
    flab/.style={font=\scriptsize\bfseries,text=red!70!black,fill=white,inner sep=1.5pt,rounded corners=1pt}
]

%------------------------------------------------
% Vertices of P_1 = R
%------------------------------------------------
\coordinate (v125) at (0,0,0);
\coordinate (v128) at (0,0,2);
\coordinate (v156) at (0,2,0);
\coordinate (v167) at (0,3,1);
\coordinate (v178) at (0,3,2);

\coordinate (v235) at (1,0,0);
\coordinate (v238) at (1,0,2);

\coordinate (v345) at (2,1,0);
\coordinate (v348) at (2,1,2);

\coordinate (v456) at (2,2,0);
\coordinate (v467) at (2,3,1);
\coordinate (v478) at (2,3,2);

%------------------------------------------------
% Filled faces (back to front)
%------------------------------------------------
\fill[facetfill] (v345)--(v348)--(v478)--(v467)--(v456)--cycle; % r=4
\fill[facetfill] (v235)--(v238)--(v348)--(v345)--cycle;         % r=3
\fill[facetfill] (v128)--(v178)--(v478)--(v348)--(v238)--cycle; % r=8
\fill[facetfill] (v167)--(v178)--(v478)--(v467)--cycle;         % r=7
\fill[facetfill] (v156)--(v167)--(v467)--(v456)--cycle;         % r=6
\fill[facetfill] (v125)--(v128)--(v238)--(v235)--cycle;         % r=2
\fill[facetfill] (v125)--(v235)--(v345)--(v456)--(v156)--cycle; % r=5
\fill[facetfill] (v125)--(v128)--(v178)--(v167)--(v156)--cycle; % r=1

%------------------------------------------------
% Hidden edges (chosen viewing convention)
%------------------------------------------------
\draw[hidedge] (v128)--(v178);
\draw[hidedge] (v178)--(v167);
\draw[hidedge] (v167)--(v156);
\draw[hidedge] (v167)--(v467);

%------------------------------------------------
% Visible edges
%------------------------------------------------
\draw[visedge] (v125)--(v128);
\draw[visedge] (v128)--(v238);
\draw[visedge] (v238)--(v348)--(v478)--(v467)--(v456)--(v345)--(v235)--(v125);
\draw[visedge] (v125)--(v156);
\draw[visedge] (v156)--(v456);
\draw[visedge] (v235)--(v238);
\draw[visedge] (v235)--(v345);
\draw[visedge] (v345)--(v348);
\draw[visedge] (v456)--(v467);
\draw[visedge] (v467)--(v478);
\draw[visedge] (v178)--(v478);

%------------------------------------------------
% Vertices
%------------------------------------------------
\foreach \P in {125,128,156,167,178,235,238,345,348,456,467,478}
  \node[vtx] at (v\P) {};

%------------------------------------------------
% Vertex labels
%------------------------------------------------
\node[vlab,below left=2pt]  at (v125) {$125$};
\node[vlab,left=3pt]        at (v128) {$128$};
\node[vlab,left=3pt]        at (v156) {$156$};

% main correction here:
\node[vlab,right=8pt]        at ($(v167)+(-0.03,0.00,0.00)$) {$167$};
\node[vlab,above=7pt]       at ($(v238)+(0.00,0.00,0.02)$) {$238$};

\node[vlab,above left=2pt]  at (v178) {$178$};
\node[vlab,below=3pt]       at (v235) {$235$};
\node[vlab,below right=2pt] at (v345) {$345$};
\node[vlab,above right=2pt] at ($(v348)+(0.10,0.02,0.02)$) {$348$};
\node[vlab,below right=2pt] at (v456) {$456$};
\node[vlab,right=4pt]       at ($(v467)+(0.06,0.00,0.00)$) {$467$};
\node[vlab,above right=2pt] at (v478) {$478$};

%------------------------------------------------
% Face labels
%------------------------------------------------
\node[flab] at (-0.28,1.45,1.02) {$1$};
\node[flab] at (0.43,-0.16,1.00) {$2$};
\node[flab] at (1.52,0.35,1.08) {$3$};
\node[flab] at (2.22,1.95,1.02) {$4$};
\node[flab] at (1.06,1.00,-0.08) {$5$};
\node[flab] at (1.13,2.52,0.32) {$6$};
\node[flab] at (1.26,2.95,1.32) {$7$};
\node[flab] at (1.08,1.58,2.12) {$8$};

\end{tikzpicture}

\caption{The polytope $\mathcal P_1$. 
A vertex label $ijk$ stands for the facet $\{i,j,k\}$, and a face label corresponds to a position index in the subword complex.}
\label{fig:P1-polytope}

\end{figure}
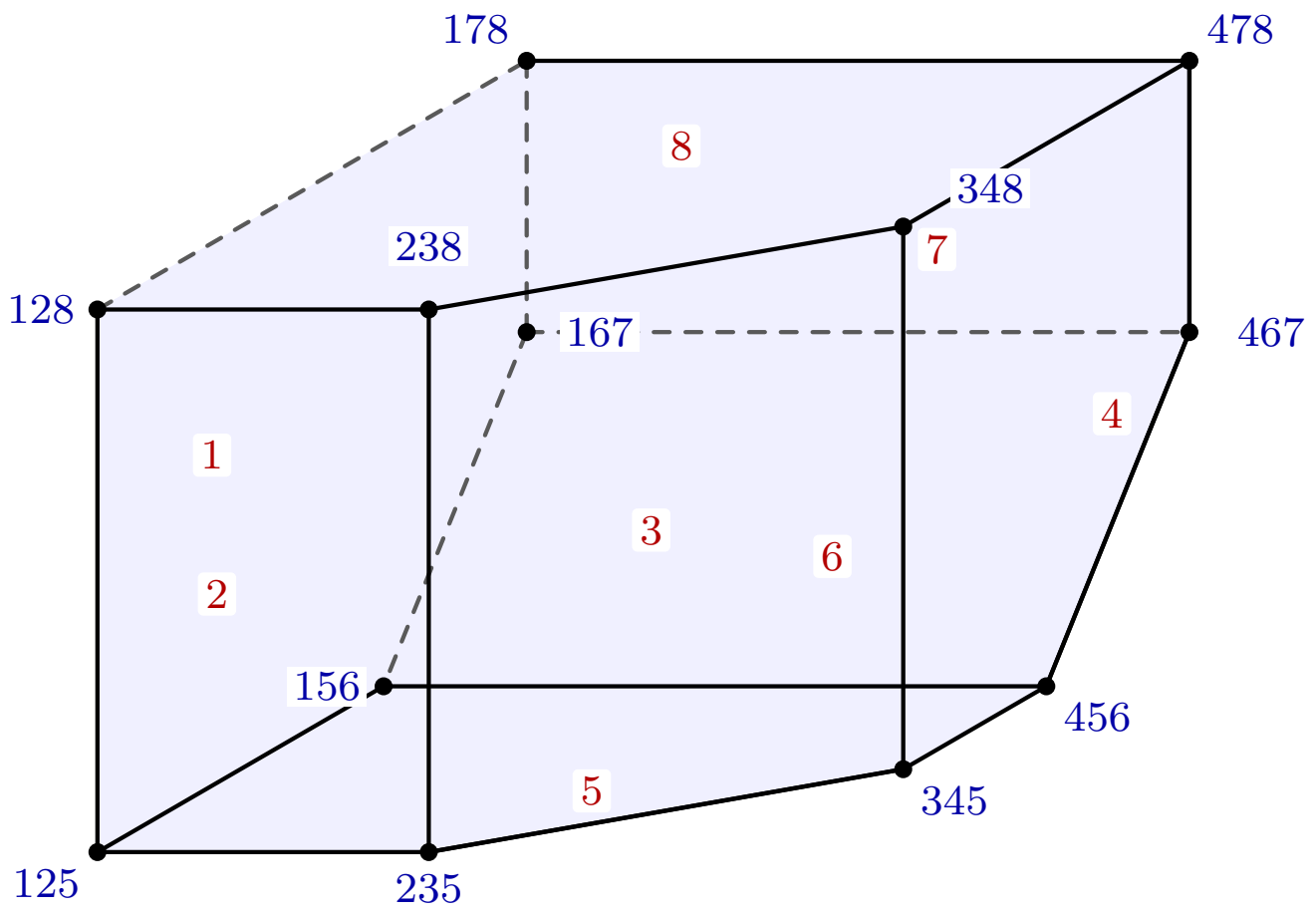

\subsection{Example 2}
\label{ex:int-small}

Now we take $W=A_3$, $Q=s_1s_2s_3s_2s_1s_3s_2s_1s_2s_3,\text{ and }
\pi=w_0=s_2s_3s_2s_1s_2s_3$.
We have $\ell(\pi)=6$, $d=4$, and $I^+=\{1,2,3,5\}$.
The root configuration $R(I^+)$ is
\[
(\alpha_1,\alpha_2,\alpha_3,\alpha_1+\alpha_2).
\]
We have
\[
M_Q=\Z A\oplus\Z U\oplus\Z V\oplus\Z W,\qquad
(A,U,V,W)=(e_1^{I^+},e_2^{I^+},e_3^{I^+},e_5^{I^+}).
\]
For coordinates $(a,u,v,w)$, the projection and its kernel are
\begin{equation}\label{eq:int-small-kernel}
\begin{aligned}
p(a,u,v,w)&=-(a+w)\alpha_1-(u+w)\alpha_2-v\alpha_3,\\
\ker p&=\Z(A+U-W).
\end{aligned}
\end{equation}

For each of the twelve index sets $K$ in Table~\ref{tab:int-base}, we use
a relabeling $\sigma:i\to i+1, i\in\{1,2,3,6,7,8\}$, $\sigma(5)=5$, $\sigma(4)=6$ and set
\[
I_{0,K}=\{1\}\cup\sigma(K),\qquad
I_{1,K}=\{10\}\cup\sigma(K).
\]
These are precisely all the facets of this subword complex. If $z_K=(u,v,w)$, their lifted
brick vectors are
\begin{equation}\label{eq:int-small-vertices}
\wt B(I_{0,K})=(0,u,v,w),\qquad
\wt B(I_{1,K})=(5+u-w,u,v,w).
\end{equation}
Thus we have all lifted brick vectors and the polytope $\mathcal P_2$ is given by
\begin{equation}\label{eq:int-small-polytope}
\mathcal P_2=\{(a,u,v,w)\mid (u,v,w)\in\mathcal P_1,\quad
0\leq a\leq5+u-w\}.
\end{equation}

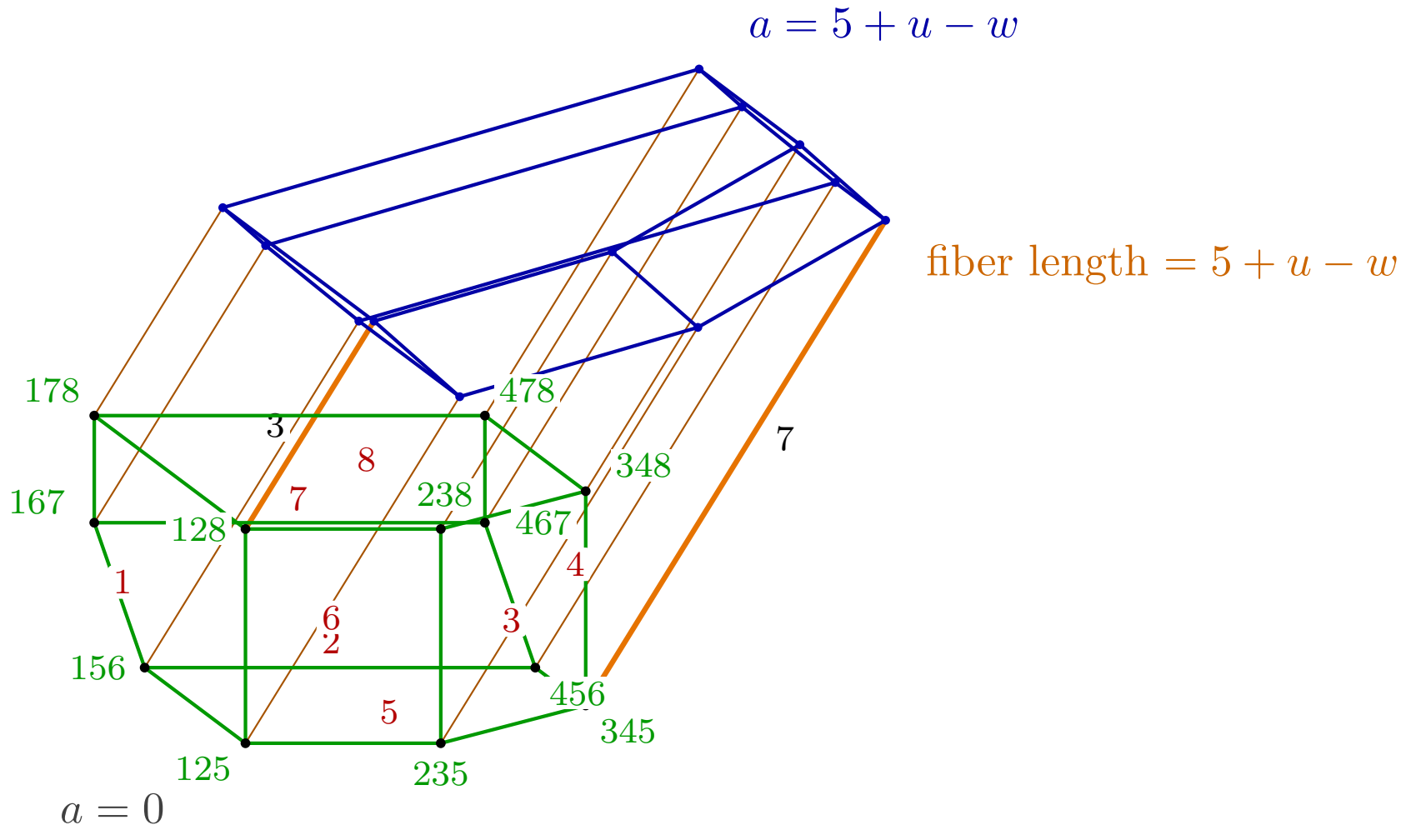
\begin{figure}[htbp]
\centering
\begin{tikzpicture}[
  x={(1.55cm,0cm)},
  y={(-0.40cm,0.30cm)},
  z={(0cm,0.85cm)},
  line cap=round,
  line join=round,
  lowedge/.style={draw=green!60!black,thick},
  topedge/.style={draw=blue!65!black,thick},
  fiber/.style={draw=orange!75!black,thick},
  fiberthin/.style={draw=orange!65!black,thin},
  important/.style={draw=orange!90!black,very thick},
  vlabel/.style={
    font=\scriptsize\bfseries,
    text=green!60!black,
    fill=white,
    inner sep=1pt
  },
  flabel/.style={
    font=\scriptsize\bfseries,
    text=red!70!black,
    fill=white,
    inner sep=1.2pt,
    rounded corners=1pt
  }
]

%================================================

%================================================
\foreach \id/\u/\v/\w in {125/0/0/0,128/0/0/2,156/0/2/0,167/0/3/1,178/0/3/2,235/1/0/0,238/1/0/2,345/2/1/0,348/2/1/2,456/2/2/0,467/2/3/1,478/2/3/2}{%
  \pgfmathsetmacro{\fiberlength}{5+\u-\w}%
  \pgfmathsetmacro{\dx}{0.34*\fiberlength}%
  \pgfmathsetmacro{\dy}{0.55*\fiberlength}%
  \coordinate (b\id) at (\u,\v,\w);%
  \coordinate (t\id) at ($(b\id)+(\dx cm,\dy cm)$);%
}

%================================================
% Fibers in the a-direction
%================================================
\foreach \id in {125,128,156,167,178,235,238,345,348,456,467,478}{%
  \draw[fiberthin] (b\id)--(t\id);%
}

% Highlight shortest and longest fibers
\draw[important] (b128)--(t128)
  node[pos=.50,left=4pt,font=\scriptsize,fill=white,inner sep=1pt] {$3$};

\draw[important] (b345)--(t345)
  node[pos=.55,right=4pt,font=\scriptsize,fill=white,inner sep=1pt] {$7$};

%================================================
% Graph of the bottom P_1
%================================================
\foreach \i/\j in {125/128,125/156,125/235,128/178,128/238,156/167,156/456,167/178,167/467,178/478,235/238,235/345,238/348,345/348,345/456,348/478,456/467,467/478}{%
  \draw[lowedge] (b\i)--(b\j);%
}

%================================================
% Graph of the top layer
%================================================
\foreach \i/\j in {125/128,125/156,125/235,128/178,128/238,156/167,156/456,167/178,167/467,178/478,235/238,235/345,238/348,345/348,345/456,348/478,456/467,467/478}{%
  \draw[topedge] (t\i)--(t\j);%
}

%================================================
% Vertex markers
%================================================
\foreach \id in {125,128,156,167,178,235,238,345,348,456,467,478}{%
  \fill[black] (b\id) circle (1.1pt);%
  \fill[blue!70!black] (t\id) circle (1.0pt);%
}

%================================================
% Bottom-layer vertex labels (all 12)
%================================================
\node[vlabel,below left=2pt]  at (b125) {$125$};
\node[vlabel,left=3pt]        at (b128) {$128$};
\node[vlabel,left=3pt]        at (b156) {$156$};
\node[vlabel,above left=1pt]  at ($(b167)+(-0.10,0.00)$) {$167$};
\node[vlabel,above left=2pt]  at (b178) {$178$};

\node[vlabel,below=3pt]       at (b235) {$235$};
\node[vlabel,above=4pt]       at ($(b238)+(0.02,0.00)$) {$238$};

\node[vlabel,below right=2pt] at (b345) {$345$};
\node[vlabel,above right=2pt] at ($(b348)+(0.08,0.00)$) {$348$};

\node[vlabel,below right=2pt] at (b456) {$456$};
\node[vlabel,right=3pt]       at ($(b467)+(0.06,0.00)$) {$467$};
\node[vlabel,above right=2pt] at (b478) {$478$};

%================================================
% Face labels for the lower P_1
%================================================
\node[flabel] at (-0.26,1.43,1.00) {$1$};
\node[flabel] at (0.40,-0.15,1.00) {$2$};
\node[flabel] at (1.45,0.34,1.03) {$3$};
\node[flabel] at (2.18,1.90,1.00) {$4$};
\node[flabel] at (1.00,1.02,-0.07) {$5$};
\node[flabel] at (1.08,2.48,0.30) {$6$};
\node[flabel] at (1.03,2.95,1.24) {$7$};
\node[flabel] at (1.02,1.55,2.10) {$8$};

%================================================
% Labels for the two layers
%================================================
\node[font=\small,text=black!75,anchor=east,fill=white,inner sep=1pt]
  at ($(b125)+(-0.60cm,-0.52cm)$) {$a=0$};

\node[font=\small,text=blue!65!black,anchor=west,fill=white,inner sep=1pt]
  at ($(t478)+(0.35cm,0.35cm)$) {$a=5+u-w$};

\node[font=\small,text=orange!80!black,anchor=west,fill=white,inner sep=1pt]
  at ($(t345)+(0.28cm,-0.35cm)$) {fiber length $=5+u-w$};

\end{tikzpicture}

\caption{
An affine projection of the four-dimensional lifted brick polytope $\mathcal P_2$.
The green part is the polytope $\mathcal P_1$. The red indices are facets of $P_1$. The blue part is a copy of $\mathcal P_1$. The orange segments represent the $a$-direction fibers.
}
\label{fig:P2-variable-fibers-refined}
\end{figure}

\subsection{Example 3}
We take $W=\widetilde A_2$,
$Q=s_1s_0s_1s_0s_2s_1s_2s_0s_1s_0s_1$ and
$\pi=\Dem(Q)=s_1s_0s_1s_2s_1s_0s_1$.
The length of $\pi$ is $7$ and $d=11-7=4$. Take
$I^+=\{1,2,5,8\}$. We have
\[
R(I^+)=(\alpha_1,\alpha_0,\alpha_1+\delta,\alpha_0+\delta),
\text{ where }\delta=\alpha_0+\alpha_1+\alpha_2.
\]
These roots are pairwise distinct up to sign, but
\[
\alpha_1-\alpha_0-(\alpha_1+\delta)+(\alpha_0+\delta)=0.
\]
This example is double root free but not root independent. 
Set
\[
M_Q=\Z A\oplus\Z B\oplus\Z C\oplus\Z D
\text{ and }
(A,B,C,D)=(e_1^{I^+},e_2^{I^+},e_5^{I^+},e_8^{I^+}).
\]
The projection from $M_Q$ to
$\Lambda_Q=\Z\alpha_0\oplus\Z\alpha_1\oplus\Z\alpha_2$ is
\begin{equation}\label{eq:int-P3-p}
p(a,b,c,d)=-a\alpha_1-b\alpha_0-c(\alpha_1+\delta)-d(\alpha_0+\delta).
\end{equation}
Its kernel is
\[
\ker p=\Z\kappa,\qquad \kappa=A-B-C+D.
\]

All $29$ facets and the coordinates of the lifted brick vectors after normalization are listed in Table 2.
\begin{table}[H]
\centering
\caption{$\SC_{\widetilde A_2}(Q,\pi)$: facets and lifted brick vectors.}
\label{tab:int-P3-vertices}
\medskip
\small
\begin{tabular}{@{}cc@{}}
\toprule
$I$ & $\wt B(I)=(a,b,c,d)$\\
\midrule
$\{1,2,5,8\}$ & $(0,0,0,0)$\\
$\{1,2,5,11\}$ & $(0,0,0,2)$\\
$\{1,2,8,9\}$ & $(0,0,2,0)$\\
$\{1,2,9,10\}$ & $(0,0,3,1)$\\
$\{1,2,10,11\}$ & $(0,0,3,2)$\\
$\{1,5,6,8\}$ & $(0,2,0,0)$\\
$\{1,5,6,11\}$ & $(0,2,0,2)$\\
$\{1,6,7,8\}$ & $(0,3,1,0)$\\
$\{1,6,7,11\}$ & $(0,3,1,2)$\\
$\{1,7,8,9\}$ & $(0,3,2,0)$\\
$\{1,7,9,10\}$ & $(0,3,3,1)$\\
$\{1,7,10,11\}$ & $(0,3,3,2)$\\
$\{2,3,5,8\}$ & $(1,0,0,0)$\\
$\{2,3,5,11\}$ & $(1,0,0,2)$\\
$\{2,3,8,9\}$ & $(1,0,2,0)$\\
$\{2,3,9,10\}$ & $(1,0,3,1)$\\
$\{2,3,10,11\}$ & $(1,0,3,2)$\\
$\{3,4,5,8\}$ & $(2,1,0,0)$\\
$\{3,4,5,11\}$ & $(2,1,0,2)$\\
$\{3,4,8,9\}$ & $(2,1,2,0)$\\
$\{3,4,9,10\}$ & $(2,1,3,1)$\\
$\{3,4,10,11\}$ & $(2,1,3,2)$\\
$\{4,5,6,8\}$ & $(2,2,0,0)$\\
$\{4,5,6,11\}$ & $(2,2,0,2)$\\
$\{4,6,7,8\}$ & $(2,3,1,0)$\\
$\{4,6,7,11\}$ & $(2,3,1,2)$\\
$\{4,7,8,9\}$ & $(2,3,2,0)$\\
$\{4,7,9,10\}$ & $(2,3,3,1)$\\
$\{4,7,10,11\}$ & $(2,3,3,2)$\\
\bottomrule
\end{tabular}
\end{table}

For the facet $I^+$, the brick vector is
\[
B(I^+)=4\omega_0+5\omega_1+2\omega_2-5\alpha_0-7\alpha_1-\alpha_2.
\]
The lifted brick polytope is
\begin{equation}\label{eq:int-P3}
\mathcal P_3=\{(a,b,c,d)\mid 0\le a\le2,\ 0\le b\le3,\ 0\le c\le3,\ 0\le d\le2,\
a-b\le1,\ b-c\le2,\ c-d\le2\}.
\end{equation}
It is a four-dimensional simple polytope with $29$ vertices and $11$ facets. Notice that the two distinct facets $I=\{1,7,8,9\}$ and $J=\{3,4,5,11\}$ satisfy $\wt B(J)-\wt B(I)=2\kappa \text{ and } B(J)=B(I).$ Thus the ordinary brick vectors do not distinguish these two facets.